\documentclass[a4paper,twoside]{article}
\usepackage{natbib}
\usepackage{amsmath,amsfonts,amssymb,amscd,graphics}
\usepackage{theorem}
\usepackage[english]{babel}
 
\theoremstyle{plain}
\newtheorem{theo}{Theorem}
\newtheorem{definition}[theo]{Definition}

\newtheorem{prop}[theo]{Proposition}
\newtheorem{lemme}[theo]{Lemma}

{\theorembodyfont{\upshape}
}

{\theorembodyfont{\upshape}
}

{\theorembodyfont{\upshape}
\newtheorem{rem}[theo]{Remark}}

\def\demo{\noindent \textsc{Proof} \\}
\def\findem{\hfill$\square$\vskip 13pt}

\def\N{\mathbb{N}}
\def\R{\mathbb{R}}

\def\E{\mathbb{E}}
\def\F{\mathcal{F}}
\def\G{\mathcal{G}}
\def\B{\mathcal{B}}
\def\P{\mathbb{P}}
\def\cvp{\xrightarrow{\P}}
\def\cvl{\xrightarrow{\mathcal{L}}}
\def\cvps{\xrightarrow{a.s.}}
\def\cvBC{\xrightarrow{BC}}
\def\cvL1{\xrightarrow{L^1}}

\title{Convergence of values in optimal stopping }

\author{Sandrine TOLDO\footnote{{\it Email address :} {\rm sandrine.toldo@math.univ-rennes1.fr}} \\
{\it IRMAR, Universit\'e Rennes 1, Campus de Beaulieu, 35042 Rennes Cedex, France}}

\date{}

\begin{document}

\maketitle

\noindent
\textbf{Abstract :} Under the hypothesis of convergence in probability of a sequence of c\`adl\`ag processes $(X^n)_n$ to a c\`adl\`ag 
process $X$, we are interested in the convergence of corresponding values in optimal stopping. We give results under hypothesis of 
inclusion of filtrations or convergence of filtrations. \\


\section{Introduction}

Let us consider a c\`adl\`ag process $X$. Let us denote by $\F^X$ its natural filtration and by $\F$ the right-continuous associated 
filtration ($\forall t, \F_t = \F^X_{t^+}$). We denote by $\mathcal{T}_L$ the set of $\F$ stopping times bounded by $L$. \\
\indent Let $\gamma : [0, + \infty [ \times \R \to \R$ a bounded continuous function. We define the value in optimal stopping of horizon $L$ 
of the process $X$ by :
$$\Gamma(L)=\underset{\tau \in \mathcal{T}_L}{\sup}\E[\gamma(\tau,X_\tau)].$$ 

\begin{rem} 
\label{GammaX}
As it is written in \citet{cvreduites}, the value of $\Gamma(L)$ only depends on the law of $X$.
\end{rem}

We are interested in the following problem : let us consider a sequence $(X^n)_n$ of processes which converges in 
probability to a limit process $X$. For all $n$, we denote by $\F^n$ the natural filtration of $X^n$ and by $\mathcal{T}^n_L$ the set of 
$\F^n$ stopping times bounded by $L$. Then, we define the values in optimal stopping $\Gamma_n(L)$ by 
$\Gamma_n(L)=\underset{\tau \in \mathcal{T}^n_L}{\sup}\E[\gamma(\tau,X^n_\tau)].$ The main aim of this paper is to give conditions under
which $(\Gamma_n(L))_n$ converges to $\Gamma(L)$. \\
\indent In his unpublished manuscript \citep{preprintAldous}, Aldous proved that if $X$ is quasi-left continuous and if there is extended 
convergence (in law) of $((X^n,\F^n))_n$ to $(X,\F)$, then $(\Gamma_n(L))_n$ converges to $\Gamma(L)$. In their paper \citep{cvreduites}, 
Lamberton and Pag\`es obtained the same result under the hypothesis of weak extended convergence of $((X^n,\F^n))_n$ to 
$(X,\F)$, quasi-left continuity of the $X^n$'s and Aldous' criterion of tightness for $(X^n)_n$. \\
\indent As a first step, we are going to prove in section \ref{liminf} that, under very weak hypothesis, holds the inequality 
$\Gamma(L) \leqslant \liminf \Gamma_n(L)$. \\
\indent Then, to prove that $(\Gamma_n(L))_n$ converges to $\Gamma(L)$, it remains to show that 
$\Gamma(L) \geqslant \limsup \Gamma_n(L)$. This inequality is more difficult and both papers \citep{preprintAldous} and \citep{cvreduites} need
weak extended convergence to prove it. Here, we prove that it happens under the hypothesis of inclusion of filtrations $\F^n \subset \F$ or
under convergence of filtrations. \\
\indent The main idea in our proof of the inequality $\Gamma(L) \geqslant \limsup \Gamma_n(L)$ is the following. We build a sequence 
$(\tau^n)$ of $(\F^n)$ stopping times bounded by $L$. Then, we extract a convergent subsequence of $(\tau^n)$ to a random variable 
$\tau$ and, at the same time, we wish to compare $\E[\gamma(\tau,X_\tau)]$ and $\Gamma(L)$. We are going to do that through two methods.\\
\indent First, we will enlarge the space of stopping times, by considering the ran\-do\-mi\-zed stopping times and the topology introduced in
\citep{BC}. Baxter and Chacon have shown that the space of randomized stopping times for a right continuous filtration with the associated 
topology is compact. We are going to use this method in section \ref{limsup_inclusion} when we have the inclusion of the filtrations 
$\F^n \subset \F$ (it means that $\forall t \in [0,T], \F^n_t \subset \F_t$).\\
\indent When we do not have the previous inclusion, we enlarge the filtration $\F$ associated to the limiting process $X$. This method is 
used, in a slightly different way, in \citep{preprintAldous} and in \citep{cvreduites}. In section \ref{limsup_cvfiltrations}, we enlarge 
(as little as possible) the limiting filtration so that the limit $\tau^*$ of a convergent subsequence of the randomized $(\F^n)$ stopping 
times associated to the $(\tau^n)_n$ is a randomized stopping time for this enlarged filtration and we use convergence of filtrations 
instead of extended convergence. \\
\indent For technical reasons, we need Aldous' criterion of tightness for the sequence $(X^n)_n$. In section \ref{condAldous}, we are going 
to show that, if $X^n \cvp X$, Aldous' criterion of tightness for $(X^n)_n$ and quasi-left continuity of the limiting process $X$ are 
e\-qui\-va\-lent.\\
\indent Finally, in section \ref{appl}, we give applications of the convergence of values in optimal stopping to discretizations and also
to financial models. \\

In what follows, we are given a probability space $(\Omega, \mathcal{A}, \P)$. We fix a positive real $T$ and also $L$ between 0 and $T$. 
Unless otherwise specified, every $\sigma$-field is supposed to be included in $\mathcal{A}$, every process will be indexed by $[0,T]$ and
taking values in $\R$ and every filtration will be indexed by $[0,T]$. $\mathbb{D}=\mathbb{D}([0,T])$ denotes the space of c\`adl\`ag
functions from $[0,T]$ to $\R$. We endow $\mathbb{D}$ with the Skorokhod topology. \\
\indent For technical background about Skorokhod topology, the reader may refer to \citep{Bill} or \citep{JS}.\\


\section{Statement of the result of convergence of the optimal values}
\label{cvvaleurs}

The main purpose of this paper is to prove the following Theorem :

\begin{theo}
\label{cvGamma}
Let us consider a c\`adl\`ag continuous in probability process X and a sequence $(X^n)_n$ of c\`adl\`ag processes. Let $\F$ be the right continuous filtration
associated to the natural filtration of X and $(\F^n)_n$ the natural filtrations of the processes $(X^n)_n$.
We assume that $X^n \cvp X$ and that one of the following assertions holds : \\
- for all n, $\F^n \subset \F$, \\
- $\F^n \xrightarrow{w} \F$. \\
Then, $\Gamma_n(L) \xrightarrow[n \to \infty]{} \Gamma(L)$.
\end{theo}

The notion of convergence of filtrations has been defined in \citep{Hoover} :

\begin{definition}
We say that $(\F^n)$ converges weakly to $\F$ if for every $A \in \F_T$, $(\E[1_A | \F^n_.])_n$ 
converges in probability to $\E[1_A | \F_.]$ for Skorokhod topology. We denote $\F^n \xrightarrow{w} \F$.
\end{definition}

The proof of Theorem \ref{cvGamma} will be given through two steps : \\
- Step 1 : we show that $\Gamma(L) \leqslant \liminf \Gamma_n(L)$, \\
- Step 2 : we show that $\Gamma(L) \geqslant \limsup \Gamma_n(L)$. \\


\section{Proof of the inequality $\Gamma(L) \leqslant \liminf \Gamma_n(L)$}
\label{liminf}

\begin{theo}
\label{thliminf} 
Let us consider a c\`adl\`ag process $X$ such that \mbox{$\P[\Delta X_L \not= 0]=0$,} its natural filtration $\F^X$, 
a sequence of c\`adl\`ag processes $(X^n)_n$ and their natural filtrations $(\F^n)_n$. We suppose that \mbox{$X^n \cvp X$.  }
Then $\Gamma(L) \leqslant \liminf \Gamma_n(L)$.
\end{theo}

\demo
\indent The proof is broken in several steps.

\begin{lemme} 
\label{T^n}
Let $\tau$ be a $\F^X$ stopping time bounded by L taking values in a discrete set $\{t_i\}_{i \in I}$ such that 
$\P[\Delta X_{t_i} \not= 0]=0$, $\forall i$. For all $i$, we consider \mbox{$A_i = \{\tau=t_i\}$.} We define 
$\tau^n$ by : 
\mbox{$\tau^n(\omega)=\min \{t_i : i \in\{j : \E[1_{A_j}|\F^n_{t_j}](\omega) > 1/2\}\}$}, $\forall \omega$.
Then, $(\tau^n)$ is a sequence of $(\mathcal{T}^n_L)$ such that $(\tau^n,X^n_{\tau^n}) \cvp (\tau,X_\tau)$.
\end{lemme}

\demo
\indent $(\tau^n)_n$ is, by definition, a sequence of $(\F^n)$-stopping times. \\

Moreover, for all $\omega$, $\tau^n(\omega) \leqslant \max \{t_i, i \in I\} \leqslant L$ because $\tau$ is bounded by $L$. So, $(\tau^n)_n$ is a sequence 
of $(\F^n)$ stopping times bounded by $L$. \\

Let us show that $\tau^n \cvp \tau$. \\

To prove that, we are going to use the convergence of $\sigma$-fields (also defined in \citep{Hoover}) :
\begin{definition}
We say that $(\mathcal{A}^n)_n$ converges to $\mathcal{A}$ and denote $\mathcal{A}^n \to \mathcal{A}$ if for every $A \in \mathcal{A}$, 
$\E[1_A | \mathcal{A}^n] \cvp 1_A$.
\end{definition}

We have $\F^n_{t_i} \to \F^X_{t_i}, \forall i$ according to the following Lemma : 
\begin{lemme}
\label{l5}
Let $(X^n)_n$ be a sequence of c\`adl\`ag processes that converges in probability to a c\`adl\`ag process X, $(\F^n)$ the natural filtrations of the 
$X^n$'s and $\F^X$ the natural filtration of X. Then, for all t, $\F^n_t \to \F^X_t$.
\end{lemme}
\demo
Take $t \in [0,T]$. Let us fix $t_1 < \ldots < t_k \leqslant t$ such that for all $i$, $\P[\Delta X_{t_i} \not=0]=0$ and let 
$f : \R^k \to \R$ be a bounded continuous. \\
$X^n \cvp X$ and for all $i=1, \ldots, k$, $\P[|\Delta X_{t_i}| \not= 0]=0$, so
$$(X^n_{t_1}, \ldots, X^n_{t_k}) \cvp (X_{t_1}, \ldots, X_{t_k}).$$
$f$ is bounded continuous so :
\begin{equation}
\label{eqcv}
f(X^n_{t_1}, \ldots, X^n_{t_k}) \xrightarrow{L^1} f(X_{t_1}, \ldots, X_{t_k}).
\end{equation}
Take $\varepsilon >0$.
\begin{eqnarray*}
\lefteqn{\P[|\E[f(X_{t_1}, \ldots, X_{t_k})|\F^n_t] - f(X_{t_1}, \ldots, X_{t_k})| \geqslant \varepsilon]} \\
& \leqslant & \P[|\E[f(X_{t_1}, \ldots, X_{t_k})|\F^n_t] - \E[f(X^n_{t_1}, \ldots, X^n_{t_k})|\F^n_t]| \geqslant \varepsilon/2]\\
& & \ \ \ \ + \P[|\E[f(X^n_{t_1}, \ldots, X^n_{t_k})|\F^n_t] - f(X_{t_1}, \ldots, X_{t_k})| \geqslant \varepsilon/2] \\
& \leqslant & \frac{4}{\varepsilon} \E[|f(X^n_{t_1}, \ldots, X^n_{t_k}) - f(X_{t_1}, \ldots, X_{t_k})|] \\
& & \ \ \ \ \text{using Markov's inequality} \\
& \xrightarrow[n \to \infty]{} & 0 \text{~~using (\ref{eqcv}).}
\end{eqnarray*}
The conclusion comes with the following characterization of the convergence of $\sigma$-fields, whose proof use exactly same arguments as 
in the proof of Lemma 3 in \citep{cvfiltration} :
\begin{lemme}
Let $Y$ be a c\`adl\`ag process, $\mathcal{A} = \sigma( \{ Y_t, t \geqslant 0\})$ and $(\mathcal{A}^n)$ a sequence of $\sigma$-fields. 
The following conditions are equivalent : \\
$i)$ $\mathcal{A}^n \to \mathcal{A}$, \\
$ii)$ $\E[f(Y_{t_1}, \ldots, Y_{t_k}) | \mathcal{A}^n] \cvp f(Y_{t_1}, \ldots, Y_{t_k})$ for every continuous bounded function $f : \R^k \to
\R$ and $t_1, \ldots, t_k$ continuity points of $Y$.
\end{lemme}
Lemma \ref{l5} is proved. 
\findem

With this Lemma, we can prove the convergence in probability of $(\tau^n)_n$ to $\tau$.\\
Let us consider a subsequence $(\tau^{\varphi(n)})_n$ of $(\tau^n)_n$. For every $i$, the convergence of the $\sigma$-fields $(\F^n_{t_i})_n$ to 
$\F^X_{t_i}$ implies \mbox{$\E[1_{A_i}|\F^{\varphi(n)}_{t_i}] \cvp 1_{A_i}$.} By successive extractions for $i \in I$ finite, 
there exists $\psi$ such that for every $i$,
$\E[1_{A_i}|\F^{\varphi \circ \psi(n)}_{t_i}] \cvps 1_{A_i}.$ For $n$ large enough, we have $\tau^{\varphi \circ \psi(n)}=\tau$
$a.s$. Then, $\tau^{\varphi \circ \psi(n)} \cvps \tau$. It follows that $\tau^n \cvp \tau$. \\

It remains to show that $X^n_{\tau^n} \cvp X_\tau$.\\
$X^n \cvp X$ so we can find a sequence $(\Lambda^n)_n$ of random time changes such that \mbox{$\sup_t |X^n_{\Lambda^n(t)} - X_t| \cvp 0$} and 
$\sup_t |\Lambda^n(t) - t| \cvp 0$. Fix $\varepsilon >0$ and $\eta >0$. We have~:
\begin{eqnarray*}
\lefteqn{\P[|X^n_{\tau^n}-X_\tau| \geqslant \eta] }\\
& \leqslant & 
\P[|X^n_{\tau^n}-X_{(\Lambda^n)^{-1}(\tau^n)}| \geqslant \eta/2] + \P[|X_{(\Lambda^n)^{-1}(\tau^n)}-X_\tau| \geqslant \eta/2]. \\
\end{eqnarray*}
There exists $n_0$ such that for every $n \geqslant n_0$, 
$\P[\sup_t |X^n_{\Lambda^n(t)} - X_t| \geqslant \eta/2] \leqslant \varepsilon$ by choice of $(\Lambda^n)_n$. 
In particular, for every $n \geqslant n_0$,
\begin{equation}
\label{Xn}
\P[|X^n_{\tau^n}-X_{(\Lambda^n)^{-1}(\tau^n)}| \geqslant \eta/2] \leqslant \varepsilon.
\end{equation}
On the other hand, for every $i \in I$ (recall that $I$ is finite), $\P[\Delta X_{t_i} \not= 0]=0$. Then, there exists $\alpha >0$ such that for every 
$i \in I$, for every $s$,
\begin{equation}
\label{Xti}
|s-t_i| \leqslant \alpha \Rightarrow \P[|X_{t_i}-X_s| \geqslant \eta/2] \leqslant \varepsilon.
\end{equation}
$\tau^n \cvp \tau$ and $\sup_t |\Lambda^n(t) - t| \cvp 0$, so $|\tau-(\Lambda^n)^{-1}(\tau^n)| \cvp 0.$
Then, there exists $n_1$ such that for every $n \geqslant n_1$, 
\begin{equation}
\label{lambda}
\P[|\tau-(\Lambda^n)^{-1}(\tau^n)| \geqslant \alpha] \leqslant \varepsilon.
\end{equation}
Then, for every $n \geqslant n_1$, 
\begin{eqnarray}
\label{eqXtau}
\lefteqn{\P[|X_{(\Lambda^n)^{-1}(\tau^n)}-X_\tau| \geqslant \eta/2]} \\ \nonumber
& = & \P[|X_{(\Lambda^n)^{-1}(\tau^n)}-X_\tau|1_{|\tau-(\Lambda^n)^{-1}(\tau^n)| \geqslant \alpha} \geqslant \eta/2] \\ \nonumber
& &  \quad + \P[|X_{(\Lambda^n)^{-1}(\tau^n)}-X_\tau|1_{|\tau-(\Lambda^n)^{-1}(\tau^n)| < \alpha} \geqslant \eta/2] \\ \nonumber
& \leqslant &
\P[2 \sup_t |X_t| 1_{|\tau-(\Lambda^n)^{-1}(\tau^n)| \geqslant \alpha} \geqslant \eta/2] + \varepsilon \text{~~using (\ref{Xti})}\\ \nonumber
& \leqslant & \P[|\tau-(\Lambda^n)^{-1}(\tau^n)| \geqslant \alpha] + \varepsilon \\ \nonumber
& \leqslant & 2 \varepsilon \text{~~using (\ref{lambda}).}
\end{eqnarray}
So, using (\ref{Xn}) and (\ref{eqXtau}), for every $n \geqslant \max(n_0,n_1)$, 
$$\P[|X^n_{\tau^n}-X_\tau| \geqslant \eta] \leqslant 3 \varepsilon.$$
\indent Finally, $(\tau^n,X^n_{\tau^n}) \cvp (\tau,X_\tau).$ \\
Lemma \ref{T^n} is proved.
\findem

With this Lemma, we can prove that Theorem \ref{thliminf} is true for stopping times that takes a finite number of values. \\

Let us consider a subdivision $\pi$ of $[0,T]$ such that no fixed time of discontinuity of $X$ belongs to $\pi$. We denote by
$\mathcal{T}^{\pi}_L$ the set of $\F$ stopping times taking values in $\pi$ and bounded by $L$. 
Then, we define :
$$\Gamma^{\pi}(L)=\underset{\tau \in \mathcal{T}^{\pi}_L}{\sup}\E[\gamma(\tau,X_\tau)].$$ 

\begin{lemme}
\label{GammaPiL}
$\Gamma^{\pi}(L) \leqslant \liminf  \Gamma_n(L).$
\end{lemme}

\demo
Fix $\varepsilon >0$. There exists a $\F^X$ stopping time $\tau$ bounded by $L$ taking values in $\pi$ such that 
$$\E[\gamma(\tau,X_\tau)] \geqslant \Gamma^{\pi}(L) - \varepsilon.$$
According to Lemma \ref{T^n}, there exists a sequence $(\tau^n)_n$ of $\F^n$ stopping times bounded by $L$ such that 
$$(\tau^n,X^n_{\tau^n}) \cvp (\tau,X_\tau).$$
$\E[\gamma(\tau^n,X^n_{\tau^n})] \to \E[\gamma(\tau,X_\tau)]$ because $\gamma$ is bounded continuous. Moreover, by definition, 
for every $n$, $\E[\gamma(\tau^n,X^n_{\tau^n})] \leqslant \Gamma_n(L)$. Il follows that  
$$\liminf \E[\gamma(\tau^n,X^n_{\tau^n})] \leqslant \liminf \Gamma_n(L).$$
But, $\liminf \E[\gamma(\tau^n,X^n_{\tau^n})] = \E[\gamma(\tau,X_\tau)] \geqslant \Gamma^{\pi}(L) - \varepsilon.$
So, $$\Gamma^{\pi}(L) - \varepsilon \leqslant \liminf \Gamma_n(L), \forall \varepsilon >0.$$
Then, $\Gamma^{\pi}(L) \leqslant \liminf \Gamma_n(L).$
\findem

It remains to link the values of optimal stopping for stopping times taking values in finite subdivisions and $\Gamma(L)$.

\begin{lemme}
\label{cvGammaPik}
Let us consider an increasing sequence $(\pi^k)_k$ of subdivisions without fixed times of continuity of $X$ such that $L \in \pi^k$ 
for every k (it is possible because $\P[\Delta X_L \not= 0]=0$) and 
$|\pi^k| \xrightarrow[k \to +\infty]{} 0$. Then $\Gamma^{\pi^k}(L) \xrightarrow[k \to +\infty]{} \Gamma(L)$.
\end{lemme}

\demo
$(\Gamma^{\pi^k}(L))_k$ is an increasing sequence bounded from above by $\Gamma(L)$. So $(\Gamma^{\pi^k}(L))_k$ converges to a limit $l$ 
with $l \leqslant \Gamma(L)$. Let us show that $l = \Gamma(L)$.\\
Fix $\varepsilon >0$. \\
We can find $\tau \in \mathcal{T}_L$ such that $$\E[\gamma(\tau,X_\tau)] \geqslant \Gamma(L) - \varepsilon.$$ 
We denote $\pi^k=\{t^k_1, \ldots , t^k_{K_k}\}$. Then, let us consider 
$$\tau^k = \sum_{i=1}^{K_k-1} t_{i+1}^k 1_{t_i^k < \tau \leqslant t_{i+1}^k}.$$
For every $k$, $\tau^k \in \mathcal{T}^{\pi^k}_L$ because $\tau$ is bounded by $L$ and $L \in \pi^k$. 
Since $|\pi^k| \to 0$, we have $\tau^k \cvp \tau$. Moreover, $\tau^k \geqslant \tau$ and $X$
is right-continuous, so $X_{\tau^k} \cvp X_\tau$. $\gamma$ is bounded continuous, so 
$$\E[\gamma(\tau^k, X_{\tau^k})] \xrightarrow[k \to \infty]{} \E[\gamma(\tau, X_\tau)].$$
But, for every $k$, $\Gamma^{\pi^k}(L) \geqslant \E[\gamma(\tau^k, X_{\tau^k})]$. It follows that  
$$l \geqslant \E[\gamma(\tau, X_\tau)] \geqslant \Gamma(L) - \varepsilon.$$
This is true for every $\varepsilon >0$, so $l \geqslant \Gamma(L)$.\\
Then, $\Gamma^{\pi^k}(L) \xrightarrow[k \to +\infty]{} \Gamma(L)$ and Lemma \ref{cvGammaPik} si proved.
\findem

At last, Theorem \ref{thliminf} follows from Lemmas \ref{GammaPiL} and \ref{cvGammaPik}.
\findem

\begin{rem}
If $\P[\Delta X_L \not= 0] > 0$, the result may not hold any longer. Let us give an example when $L = 1/2$. We consider some processes $x$ 
and $(x^n)$ defined on $[0,1]$ by $x_t = 1_{[1/2, 1]}(t)$ and $x^n_t = 1_{[1/2+1/n, 1]}(t)$, $\forall t$. Let us consider 
$\gamma : [0, + \infty[ \times \R \to \R$ such that $\gamma(t,y)=y \wedge 2$. $\gamma$ is a continuous bounded function. We want to compare
$\Gamma(1/2)$ and the limit of $\Gamma_n(1/2)$ when $n$ goes to $+ \infty$.\\
We have : $\Gamma(1/2) = \underset{\tau \in \mathcal{T}_{1/2}}{\sup}\E[\gamma(\tau,x_\tau)] = \underset{t \leqslant 1/2}{\sup} x_t = 1$. \\
On the other hand, for every $n$, $\Gamma_n(1/2) = \underset{t \leqslant 1/2}{\sup} x^n_t = 0$. \\
So $\liminf \Gamma_n(1/2) =0 < 1 = \Gamma(1/2)$.
\end{rem}

\begin{rem}
\label{rqF}
The Theorem remains true if we replace $\F^X$ by the right continuous filtration associated to $\F$ ($\forall t, \F_t=\F^X_{t+}$) and if we
take the $\Gamma(L)$ associated to $\F$.
\end{rem}


\section{Aldous' criterion for tightness}
\label{condAldous}

In his papers \citep{Aldous_Stop1} and \citep{Aldous_Stop2}, Aldous deals with the following criterion for tightness :
\begin{equation}
\label{Aldous}
\forall \varepsilon > 0, \underset{\delta \downarrow 0}{\lim} \ \underset{n \to +\infty}{\limsup} \ 
        \underset{S, T \in \mathcal{T}_L^n, S \leqslant T \leqslant S+\delta}{\sup} \ 
                \P[| X^n_{S} - X^n_{T} | \geqslant \varepsilon] = 0.
\end{equation}
He gives many results which links that criterion and weak convergence of sequences of processes. \\

In his unpublished manuscript \citep{preprintAldous}, Aldous shows the following result (Corollary 16.23) which links convergence of stopping 
times to convergence of processes :

\begin{prop}
\label{cvlXnta}
Let us consider a sequence of c\`adl\`ag processes $(X^n)_n$ that converges in law to a c\`adl\`ag process X. We denote by $\F^n$ the natural filtrations
of the processes $X^n$ and by $\F$ the right continuous natural filtration of the process X. 
Let us consider a sequence $(\tau^n)_n$ of $(\F^n)$-stopping times that converges in law to a random variable $V$.
We suppose that we have the join convergence in law of $((\tau^n,X^n))_n$ to $(V,X)$ and that Aldous' criterion for tightness (\ref{Aldous}) 
is filled. Then $(\tau^n,X^n_{\tau^n}) \cvl (V,X_V)$.
\end{prop}

\demo
We just give the sketch of Aldous' proof. \\
\indent If $\P[\Delta X_V \not= 0]=0$, using the Skorokhod representation Theorem, we can prove that $(\tau^n,X^n_{\tau^n}) \cvl (V,X_V)$. \\
\indent If $\P[\Delta X_V \not= 0]\not=0$, we can find a decreasing sequence $(\delta_k)_k$ of reals that converges to 0 and such that for 
every $k$, $\P[\Delta X_{V+\delta_k} \not= 0]=0$. \\
Let us take $f:\R^2 \to \R$ bounded and continuous.
\begin{eqnarray*}
\lefteqn{|\E[f(\tau^{n},X^n_{\tau^{n}}) - f(V,X_V)]|} \\
& \leqslant & |\E[f(\tau^{n},X^n_{\tau^{n}}) - f(\tau^{n}+\delta_k,X^n_{\tau^{n}+\delta_k})]| \\
& & \quad + |\E[f(\tau^{n}+\delta_k,X^n_{\tau^{n}+\delta_k})-f(V+\delta_k,X_{V+\delta_k})]| \\
& & \quad \quad + |\E[f(V+\delta_k,X_{V+\delta_k})- f(V,X_V)]|.
\end{eqnarray*}
But : \\
- $\forall k$, $\limsup_{n \to + \infty} \E[f(\tau^{n}+\delta_k,X^n_{\tau^{n}+\delta_k})-f(V+\delta_k,X_{V+\delta_k})]=0$ 
because $\P[\Delta X_{V+\delta_k} \not= 0]=0$,  \\
- $\lim_{k \to +\infty} \E[f(V+\delta_k,X_{V+\delta_k})- f(V,X_V)]=0$ because $X$ is right-continuous, \\
- $\lim_{k \to +\infty} \limsup_{n \to + \infty} \E[f(\tau^{*,n},X^n_{\tau^{*,n}}) - f(\tau^{*,n}+\delta_k,X^n_{\tau^{*,n}+\delta_k})]=0$ 
using Aldous' criterion. \\
The result follows.
\findem

\begin{rem}
We will see in Proposition \ref{cvlXntar} an analogous result in the case of randomized stopping times.
\end{rem}

The following caracterization of Aldous' Criterion is probably widely known, however I do not know of any reference to a proof of it, so I
give one of my own here.

\begin{prop}
\label{hypA}
Let us consider a sequence of c\`adl\`ag processes $(X^n)_n$ and a c\`adl\`ag process X such that $X^n \cvp X$. 
The following conditions are equivalent~:\\
i) X is continuous in probability everywhere, ie for every t $\P[\Delta X_t \not= 0]=0$, \\ 
ii) Aldous' criterion for tightness (\ref{Aldous}) is filled. 
\end{prop}

\demo
$i) \Rightarrow ii).$ 
Let $\delta > 0$. Let $(T^n)_n$ and $(S_n)_n$ be two sequences of $\mathcal{T}_L^{n}$ such that 
for every $n$, \mbox{$S^n \leqslant T^n \leqslant S^n + \delta$.} Let $\varepsilon >0$ and $\eta > 0$. \\
$X^n \cvp X$ so we can find a sequence of random time changes $(\Lambda^n)_n$ such that \mbox{$\sup_t |X^n_{\Lambda^n (t)} -X_t| \cvp 0$.} 
Then there exists $n_0$ such that 
$$\forall n \geqslant n_0, \P[\sup_t |X^n_{\Lambda^n (t)} -X_t| \geqslant \eta/3] \leqslant \varepsilon.$$
We have :
\begin{eqnarray}
\label{equ1}
\lefteqn{\P[| X^n_{S^n} - X^n_{T^n} | \geqslant \eta]} \nonumber \\
& \leqslant & \P[| X^n_{S^n} - X_{(\Lambda^n)^{-1}(S^n)} | \geqslant \eta/3] \nonumber \\
& & \quad + \P[|X_{(\Lambda^n)^{-1}(S^n)} - X_{(\Lambda^n)^{-1}(T^n)}| \geqslant \eta/3] \nonumber \\
& & \quad \quad + \P[|X_{(\Lambda^n)^{-1}(T^n)} - X^n_{T^n} | \geqslant \eta/3]
\end{eqnarray}
But, for every $n \geqslant n_0$,
\begin{equation}
\label{equ2}
\P[| X^n_{S^n} - X_{(\Lambda^n)^{-1}(S^n)} | \geqslant \eta/3]
\leqslant \P[\sup_t |X^n_{\Lambda^n (t)} -X_t| \geqslant \eta/3] \leqslant \varepsilon,
\end{equation}
and
\begin{equation}
\label{equ3}
\P[|X_{(\Lambda^n)^{-1}(T^n)} - X^n_{T^n} | \geqslant \eta/3]
\leqslant \P[\sup_t |X^n_{\Lambda^n (t)} -X_t| \geqslant \eta/3] \leqslant \varepsilon.
\end{equation}
It remains to show that :
$$\lim_{\delta \downarrow 0} \limsup_{n \to +\infty}
\P[|X_{(\Lambda^n)^{-1}(S^n)} - X_{(\Lambda^n)^{-1}(T^n)}| \geqslant \eta/3]=0.$$
$X$ is a c\`adl\`ag process, so there exists $\theta > 0$ such that 
$$\P[w'(X,\theta) \geqslant \eta/12] \leqslant \varepsilon,$$
where $\forall x \in \mathbb{D}, w'(x,\delta)=\inf_{\{t_i\} \in F_\delta} \max_{1 \leqslant i \leqslant v} w(x,[t_{i-1}, t_i[)$,
$F_\delta$ is the set of subdivisions $\{t_i\}_{1 \leqslant i \leqslant v}$ of $[0,T]$ such that $\forall i, t_i-t_{i-1} > \delta$ and $w$ is
the modulus of continuity $w(x,[t_{i-1}, t_i[)=\sup\{|x_t-x_s|, t_{i-1}<s<t<t_i\}$ \citep[see e.g.][Section 12]{Bill}.\\
By defintion of $w'$, there exists a subdivision $\{t_k\}$ such that 
$$\forall k, |t_{k+1}-t_k| \geqslant \theta \text{~~and~~} \P[\max_k w(X,[t_k,t_{k+1}[) \geqslant \eta/12] \leqslant 2 \varepsilon.$$
On the other hand,
\begin{eqnarray*}
\lefteqn{\P[|(\Lambda^n)^{-1}(T^n) - (\Lambda^n)^{-1}(T^n+\delta)| \geqslant \theta]} \\
& \leqslant & \P[|(\Lambda^n)^{-1}(T^n) - T^n| \geqslant \theta/3] + \P[|T^n - (T^n + \delta)| \geqslant \theta/3] \\
& & \quad + \P[|T^n + \delta - (\Lambda^n)^{-1}(T^n + \delta)| \geqslant \theta/3] \\
& \leqslant & 2 \P[\sup_t |(\Lambda^n)^{-1}(t)-t| \geqslant \theta/3]  
\text{~~for every $\delta < \theta/3$.}
\end{eqnarray*}
$\sup_t |(\Lambda^n)^{-1}(t)-t| \cvp 0$ , so there exists $n_1$ such that 
$$\forall n \geqslant n_1, \P[\sup_t |(\Lambda^n)^{-1}(t)-t| \geqslant \theta/3] \leqslant \varepsilon .$$
Then, for every $n \geqslant n_1$, for every $\delta < \theta/3$,
\begin{equation}
\label{equ*}
\P[|(\Lambda^n)^{-1}(T^n) - (\Lambda^n)^{-1}(S^n)| \geqslant \theta] \leqslant 3\varepsilon .
\end{equation}
So, for every $n \geqslant n_1$, for every $\delta < \theta/3$, 
\begin{eqnarray*}
\lefteqn{\P[|X_{(\Lambda^n)^{-1}(S^n)} - X_{(\Lambda^n)^{-1}(T^n)}| \geqslant \eta/3]} \\
& = & \P[|X_{(\Lambda^n)^{-1}(S^n)} - X_{(\Lambda^n)^{-1}(T^n)}|
1_{|(\Lambda^n)^{-1}(T^n) - (\Lambda^n)^{-1}(S^n)| < \theta} \geqslant \eta/3] \\
& & \quad + \P[|X_{(\Lambda^n)^{-1}(S^n)} - X_{(\Lambda^n)^{-1}(T^n)}|
1_{|(\Lambda^n)^{-1}(T^n) - (\Lambda^n)^{-1}(S^n)| \geqslant \theta} \geqslant \eta/3]. 
\end{eqnarray*}
But,
\begin{eqnarray*}
\lefteqn{\P[|X_{(\Lambda^n)^{-1}(S^n)} - X_{(\Lambda^n)^{-1}(T^n)}|
1_{|(\Lambda^n)^{-1}(T^n) - (\Lambda^n)^{-1}(S^n)| < \theta} \geqslant \eta/3]}\\
& \leqslant & \P[(2 \max_k w(X,[t_k,t_{k+1}[) + \max_k |\Delta X_{t_k}|) \geqslant \eta/3] \\
& \leqslant & \P[\max_k w(X,[t_k,t_{k+1}[) \geqslant \eta/12] + \P[\max_k |\Delta X_{t_k}|\geqslant \eta/6]\\
& \leqslant & 2 \varepsilon + \sum_k \P[|\Delta X_{t_k}|\geqslant \eta/6]\\
& \leqslant & 2 \varepsilon \text{~~because $X$ has no fixed time of discontinuity}
\end{eqnarray*}
and 
\begin{eqnarray*}
\lefteqn{\P[|X_{(\Lambda^n)^{-1}(S^n)} - X_{(\Lambda^n)^{-1}(T^n)}|
1_{|(\Lambda^n)^{-1}(T^n) - (\Lambda^n)^{-1}(S^n)| \geqslant \theta} \geqslant \eta/3]} \\
& \leqslant & \P[ 2 \sup_{t}|X_t| 1_{|(\Lambda^n)^{-1}(T^n) - (\Lambda^n)^{-1}(S^n)| \geqslant \theta} \geqslant \eta/3]
\\
& \leqslant & \P[|(\Lambda^n)^{-1}(T^n) - (\Lambda^n)^{-1}(S^n)| \geqslant \theta] \\
& \leqslant & 3 \varepsilon \text{~~using (\ref{equ*}).}
\end{eqnarray*}
So for every $n \geqslant n_1$, for every $\delta < \theta/3$,
\begin{equation}
\label{equ4}
\P[|X_{(\Lambda^n)^{-1}(S^n)} - X_{(\Lambda^n)^{-1}(T^n)}| \geqslant \eta/3] \leqslant 5\varepsilon .
\end{equation}
Finally, using inequalities (\ref{equ1}), (\ref{equ2}), (\ref{equ3}) and (\ref{equ4}), 
for every $n \geqslant \max(n_0, n_1)$, for every $\delta < \theta/3$, 
$$\P[| X^n_{S^n} - X^n_{T^n} | \geqslant \eta] \leqslant 7\varepsilon .$$ 
$n_0$, $n_1$ and $\theta$ do not depend on $(T_n)_n$ and $(S_n)_n$. 
Then, for every $n \geqslant \max(n_0, n_1)$, for every $\delta < \theta/3$, 
$$\sup_{S, T \in \mathcal{T}_L^n, S \leqslant T \leqslant S+\delta} \P[| X^n_{S} - X^n_{T} | \geqslant \eta] \leqslant 7\varepsilon .$$
Aldous' criterion follows.\\

\medskip 

\noindent $ii) \Rightarrow i).$ 
We suppose that there exists $t_0$ such that $\P[\Delta X_{t_0} \not= 0] >0.$ \\
Let $\varepsilon >0$ and $\eta >0$ be such that $\P[|\Delta X_{t_0}|  \geqslant 2\varepsilon] \geqslant 2\eta.$ \\
$X^n \cvp X$ so we can find a random sequence $(t^n)_n$ such that $t^n \cvp t_0$ and \mbox{$\Delta X^n_{t^n} \cvp \Delta X_{t_0}$.} There exists
$n_0$ such that for every $n \geqslant n_0$, 
\begin{equation}
\label{eq-1delta}
\P[|t^n-t_0| \geqslant \delta/2] \leqslant \eta/2 \text{~~and~~} 
\P[|\Delta X^n_{t^n} - \Delta X_{t_0}| \geqslant \varepsilon] \leqslant \eta/2.
\end{equation}
We are going to show that for every $n \geqslant n_0$, for $\delta$ large enough,
$$\P[|X^n_{t_0 + \delta/2} - X^n_{t_0 - \delta/2}| \geqslant \varepsilon/3] \geqslant \eta/2.$$
Then, for every $n \geqslant n_0$, 
\begin{eqnarray}
\label{eqdelta}
\lefteqn{\P[|\Delta X^n_{t^n}| \geqslant \varepsilon] } \nonumber \\
& = & \P[|\Delta X^n_{t^n}|1_{|t^n-t_0| \geqslant \delta/2} \geqslant \varepsilon] + 
\P[|\Delta X^n_{t^n}|1_{|t^n-t_0| < \delta/2} \geqslant \varepsilon]  \\
& \leqslant & \P[|t^n-t_0| \geqslant \delta/2] + \P[|X^n_{t^n} - X^n_{t_0 + \delta/2}|1_{|t^n-t_0| < \delta/2} \geqslant \varepsilon/3] 
\nonumber \\
&&+ \P[|X^n_{t_0 + \delta/2} - X^n_{t_0 - \delta/2}| \geqslant \varepsilon/3] 
+ \P[|X^n_{t_0 - \delta/2} - X^n_{t^{n}-}|1_{|t^n-t_0| < \delta/2} \geqslant \varepsilon/3] \nonumber
\end{eqnarray}
$(X^n)_n$ is tight. So, we can find $\delta_0>0$ and $n_1 \in \N$ such that for every $\delta \leqslant \delta_0$, for every $n \geqslant
n_1$, 
$$\P[w'(X^n,\delta) \geqslant \varepsilon/6] \leqslant \eta/6.$$
Then, we can find a finite subdivision $\{t_k\}$ such that 
$$\forall k, t_{k+1}-t_k \geqslant \delta \text{~~and~~} \P[\max_k w(X^n,[t_k,t_{k+1}[) \geqslant \varepsilon/3] \leqslant \eta/4.$$
We know that for every $n \geqslant \max(n_0,n_1)$,
\begin{eqnarray*}
2\eta & \leqslant & \P[|\Delta X_{t_0}|  \geqslant 2\varepsilon] \\
& \leqslant & \P[|\Delta X^n_{t^n} - \Delta X_{t_0}| \geqslant \varepsilon] + \P[|\Delta X^n_{t^n}| \geqslant \varepsilon] \\
& \leqslant & \eta/2 + \P[|\Delta X^n_{t^n}| \geqslant \varepsilon].
\end{eqnarray*}
In particular, for every $n \geqslant \max(n_0,n_1)$,
\begin{equation}
\label{eq0delta}
\P[|\Delta X^n_{t^n}| \geqslant \varepsilon] \geqslant 3\eta/2.
\end{equation}
So, for every $n \geqslant \max(n_0,n_1)$, $t^n \in \{t_k\}$. \\
Then, for every $\delta \leqslant \delta_0$, for every $n \geqslant \max(n_0,n_1)$,
\begin{equation}
\label{eq1delta}
\P[|X^n_{t^n} - X^n_{t_0 + \delta/2}|1_{|t^n-t_0| < \delta/2} \geqslant \varepsilon/3]
\leqslant \P[\max_k w(X^n,[t_k,t_{k+1}[) \geqslant \varepsilon/3] \leqslant \eta/4.
\end{equation}
On the same way,
\begin{equation}
\label{eq2delta}
\P[|X^n_{t_0 - \delta/2} - X^n_{t^{n}-}|1_{|t^n-t_0| < \delta/2} \geqslant \varepsilon/3] \leqslant \eta/4.
\end{equation}
Finally, using (\ref{eqdelta}) and inequalities (\ref{eq-1delta}), (\ref{eq0delta}), (\ref{eq1delta}) and (\ref{eq2delta}), 
for every $\delta \leqslant \delta_0$, for every $n \geqslant \max(n_0,n_1)$,
$$3\eta/2 \leqslant \P[|\Delta X^n_{t^n}| \geqslant \varepsilon] \leqslant \eta/2 + \eta/4 + 
\P[|X^n_{t_0 + \delta/2} - X^n_{t_0 - \delta/2}| \geqslant \varepsilon/3] + \eta/4.$$
So, for every $\delta \leqslant \delta_0$, for every $n \geqslant \max(n_0,n_1)$,
\begin{eqnarray*}
\eta/2 & \leqslant & \P[|X^n_{t_0 + \delta/2} - X^n_{t_0 - \delta/2}| \geqslant \varepsilon/3] \\
 & \leqslant & \underset{S, T \in \mathcal{T}_L^n, S \leqslant T \leqslant S+\delta}{\sup} \ 
                \P[| X^n_{T^n+\delta} - X^n_{T^n} | \geqslant \varepsilon/3].
\end{eqnarray*}
Taking the $\limsup$ when $n$ tends to infinity and the limit when $\delta$ decreases to 0, we have :
$$\eta/2 \leqslant \underset{\delta \downarrow 0}{\lim} \ \underset{n \to +\infty}{\limsup} \ 
        \underset{S, T \in \mathcal{T}_L^n, S \leqslant T \leqslant S+\delta}{\sup} \ 
                \P[| X^n_{T^n+\delta} - X^n_{T^n} | \geqslant \varepsilon/3],$$
which is in contradiction with Aldous' criterion.
The result follows.
\findem


\section{Proof of the inequality $\Gamma(L) \geqslant \limsup \Gamma_n(L)$ when for every $n$, $\F^n \subset \F$}
\label{limsup_inclusion}

\subsection{Randomized stopping times }
\label{tar}

The notion of randomized stopping times has been introduced in \citep{BC} and this notion has been used in \citep{Meyer}
 under the french name "temps d'arr\^et flous". \\

We are given a filtration $\F$. Let us denote by $\mathcal{B}$ the Borel $\sigma$-field on $[0,1]$. Then, we define the filtration 
$\mathcal{G}$ on $\Omega \times [0,1]$ such that $\forall t$, $\mathcal{G}_t=\F_t \times \mathcal{B}$. A map 
\mbox{$\tau : \Omega \times [0,1] \to [0,+\infty]$}
is called a randomized $\F$ stopping time if $\tau$ is a $\mathcal{G}$ stopping time. We denote by $\mathcal{T}^*$ the set of randomized 
stopping times and by $\mathcal{T}^*_L$ the set of randomized stopping times bounded by $L$. $\mathcal{T}$ is included in $\mathcal{T}^*$ 
and the application $\tau \mapsto \tau^*$, where $\tau^*(\omega,t)=\tau(\omega)$ for every $\omega$ and every $t$, maps $\mathcal{T}$ into 
$\mathcal{T}^*$. In the same way, $\mathcal{T}_L$ is included in $\mathcal{T}^*_L$. \\
\indent On the space $\Omega \times [0,1]$, we put the probability measure $\P \otimes \mu$ where $\mu$ is Lebesgue's measure on
$[0,1]$. In their paper \citep{BC}, Baxter and Chacon define the convergence of randomized stopping times by the following :
$$\tau^{*,n} \xrightarrow{BC} \tau^* \text{~~iff~~} \forall f \in \mathcal{C}_b ([0,\infty]), \forall Y \in L^1(\Omega, \F, \P), 
\E[Yf(\tau^{*,n})] \to \E[Yf(\tau^*)],$$
where $\mathcal{C}_b ([0,\infty])$ is the set of bounded continuous functions on $[0,\infty]$.\\
\indent Taking $Y=1$, we note that this convergence implies the "usual" convergence in law.\\
\indent This notion is a particular case of "stable convergence" introduced in \citep{Renyi} and studied in \citep{JacodMemin}. 
This is the link between convergence in probability and stable convergence that we are going to use~:
\begin{lemme}
\label{cvpcvBC}
Let us consider a sequence $(\tau^n)_n$ of $\F$ stopping times that converges in probability to $\tau$. Then the sequence $(\tau^{*,n})_n$ 
where \mbox{$\tau^{*,n}(\omega,t)=\tau^n(\omega)$} $\forall \omega$, $\forall t$, converges in Baxter and Chacon's way to $\tau^*$ where  
$\tau^{*}(\omega,t)=\tau(\omega)$, $\forall \omega$, $\forall t$.
\end{lemme}

\indent One of the main interests of this notion is, as it is shown in \citep[Theorem 1.5]{BC}, that the set of
randomized stopping times for a right continuous filtration is compact for Baxter and Chacon's topology. \\

The following Proposition is the main argument in the proof of Theorem \ref{thlimsup1} below.

\begin{prop}
\label{sssuitetar}
Let us consider a sequence of filtrations $(\F^n)$ and a right continuous filtration $\F$ such that $\forall n$, $\F^n \subset \F$. 
Let $(\tau^n)_n$ be a sequence of $(\mathcal{T}_L^{n})_n$. Then, there exists a randomized $\F$ stopping time $\tau^*$ and a subsequence
$(\tau^{\varphi(n)})_n$ such that $\tau^{*,\varphi(n)} \xrightarrow{BC} \tau^*$ where for every $n$, $\tau^{*,n}(\omega,t)=\tau^n(\omega)$ 
$\forall \omega$, $\forall t$.
\end{prop}

\demo
For every $n$, $\F^n \subset \F$, $(\tau^n)_n$ is a sequence of $\F$ stopping times so, by definition, $(\tau^{*,n} )_n$ is a sequence of
randomized $\F$ stopping times. According to \citep[Theorem 1.5]{BC}, we can find a randomized $\F$
stopping time $\tau^*$ and a subsequence $(\tau^{\varphi(n)})_n$ such that $\tau^{*,\varphi(n)} \xrightarrow{BC} \tau^*$.
\findem

\indent Now, we define $X_{\tau^*}$ by 
$X_{\tau^*}(\omega,v)=X_{\tau^*(\omega,v)}(\omega)$, for every $(\omega,v) \in \Omega \times [0,1]$. Then, we can prove the
following Lemma :
\begin{lemme}
\label{GammaGamma*}
Let us consider $\Gamma^*(L)=\underset{\tau^* \in \mathcal{T}^*_L}{\sup}\E[\gamma(\tau^*,X_{\tau^*})]$. Then $\Gamma^*(L)=\Gamma(L)$.
\end{lemme}

\demo
- $\mathcal{T}_L$ is included $\mathcal{T}^*_L$. Then $\Gamma(L) \leqslant \Gamma^*(L)$. \\
- Let $\tau^* \in \mathcal{T}^*_L$. We consider, for every $v$, $\tau_v(\omega)=\tau^*(\omega,v), \forall \omega$. \\
For every $v \in [0,1]$, for every $t \in [0,T]$,
$$\{\omega : \tau_v(\omega) \leqslant t\} \times \{v\} = \{ (\omega, x) : \tau^*(\omega,x) \leqslant t\} \cap (\Omega \times \{v\}).$$
But, $\{ (\omega, x) : \tau^*(\omega,x) \leqslant t\} \in \F_t \times \mathcal{B}$ because $\tau^*$ is a randomized $\F$ stopping time and  
$\Omega \times \{v\} \in \F_t \times \mathcal{B}$. So, $\{\omega : \tau_v(\omega) \leqslant t\} \times \{v\} \in \F_t \times
\mathcal{B}$. Consequently, 
$$\{\omega : \tau_v(\omega) \leqslant t\} \in \F_t.$$
Then, for every $v$, $\tau_v$ is a $\F$ stopping time bounded by $L$. We have :
\begin{eqnarray*}
\E[\gamma(\tau^*,X_{\tau^*})] & = & \int_\Omega \int_0^1 \gamma(\tau^*(\omega,v),X_{\tau^*(\omega,v)}(\omega))d\P(\omega)dv \\
& = & \int_0^1 \left( \int_\Omega \gamma(\tau^*(\omega,v),X_{\tau_v(\omega)}(\omega))d\P(\omega) \right) dv \\
& = & \int_0^1 \E[\gamma(\tau_v,X_{\tau_v})] dv \\
& \leqslant & \Gamma(L) \text{~~because, for every $v$, $\tau_v \in \mathcal{T}_L$.}
\end{eqnarray*}
Taking the $\sup$ for $\tau^*$ in $\mathcal{T}^*_L$, we get $\Gamma^*(L) \leqslant \Gamma(L)$. \\
Lemma \ref{GammaGamma*} is proved.
\findem

We have an analogous of Proposition \ref{cvlXnta} in the setting of randomized stopping times.

\begin{prop}
\label{cvlXntar}
Let us consider a sequence $(X^n)_n$ of c\`adl\`ag processes that converges in law to a c\`adl\`ag process X, $\F^n$ the natural filtrations 
of the $X^n$'s and $\F$ the right continuous natural filtration of the process X. Let $(\tau^n)_n$ be a sequence of $(\F^n)$ stopping times
such that the associated sequence $(\tau^{*,n})_n$ of randomized stopping times ($\tau^{*,n}(\omega,t)=\tau^n(\omega)$ $\forall \omega$, 
$\forall t$) converges in law to a random variable $V$. We suppose that $(\tau^{*,n},X^n) \cvl (V,X)$ and that Aldous' criterion \ref{Aldous}
is filled. Then $(\tau^{*,n},X^n_{\tau^{*,n}}) \cvl (V,X_V)$.
\end{prop}

\demo
The proof of Proposition \ref{cvlXntar} follows the lines of the proof of \citep[Corollary 16.23]{preprintAldous} (Proposition \ref{cvlXnta}
in this paper). 
\findem

\begin{rem}
We point out that, in this Proposition, Aldous' Criterion is filled by the original -not randomized- stopping times.
\end{rem}

When $X^n \cvp X$ and when $(\tau^{*,n})_n$ is a sequence of randomized stopping times converging in the sense of Baxter and Chacon to a 
random variable $\tau^*$, we have the join convergence in law of $((X^n, \tau^{*,n}))_n$ to $(X, \tau^*)$ : 

\begin{prop}
\label{cvcouple}
Let us consider a sequence $(X^n)_n$ of c\`adl\`ag processes converging in probability to a c\`adl\`ag process X, $\F^n$ the natural 
filtrations of the $X^n$'s and $\F$ the right continuous natural filtration of the process X. Let $(\tau^{*,n})_n$ be a sequence of 
ran\-do\-mi\-zed $(\F^n)$ stopping times converging to the ran\-do\-mi\-zed stopping time $\tau$ under Baxter and Chacon's topology. \\
Then $(X^n, \tau^{*,n}) \cvl (X, \tau^*)$.
\end{prop}

\demo
- As $(X^n)_n$ and $(\tau^{*,n})_n$ are tight, $((X^n, \tau^{*,n}))_n$ is tight. \\
- We are now going to identify the limit thanks to the finite-dimensional convergence. \\
Let $k \in \N$ and $t_1 < \ldots < t_k$ such that for every $i$, $\P[\Delta X_{t_i} \not= 0] = 0$. 
Let us show that $(X^n_{t_1}, \ldots , X^n_{t_k},\tau^{*,n}) \cvl (X_{t_1}, \ldots , X_{t_k}, \tau^*)$. \\
In a first time, let us consider $f : \R^k \to \R$ and $g : \R \to \R$ bounded continuous. 
\begin{eqnarray*}
\lefteqn{ |\E[f(X^n_{t_1}, \ldots , X^n_{t_k})g(\tau^n)] - \E[f(X_{t_1}, \ldots , X_{t_k})g(\tau^*)]| } \\
& \leqslant & |\E[(f(X^n_{t_1}, \ldots , X^n_{t_k})-f(X_{t_1}, \ldots , X_{t_k}))g(\tau^{*,n})]| \\
& & \quad + |\E[f(X_{t_1}, \ldots , X_{t_k})g(\tau^{*,n})]-\E[f(X_{t_1}, \ldots , X_{t_k})g(\tau^*)]| \\
& \leqslant & \|g\|_{\infty} \E[|f(X^n_{t_1}, \ldots , X^n_{t_k})-f(X_{t_1}, \ldots , X_{t_k})|] \\
& & \quad + |\E[f(X_{t_1}, \ldots , X_{t_k})g(\tau^{*,n})]-\E[f(X_{t_1}, \ldots , X_{t_k})g(\tau^*)]| \\
\end{eqnarray*}
But, $X^n \cvp X$ and for every $i$, $\P[\Delta X_{t_i} \not= 0] = 0$ so 
\mbox{$(X^n_{t_1}, \ldots , X^n_{t_k}) \cvp (X_{t_1}, \ldots , X_{t_k})$.}
Moreover, $f$ is bounded continuous, so 
$$\E[|f(X^n_{t_1}, \ldots , X^n_{t_k})-f(X_{t_1}, \ldots , X_{t_k})|] \xrightarrow[n \to +\infty]{} 0.$$
On the other hand, by definition of Baxter and Chacon's convergence, 
$$\E[f(X_{t_1}, \ldots , X_{t_k})g(\tau^{*,n})]-\E[f(X_{t_1}, \ldots , X_{t_k})g(\tau^*)] \xrightarrow[n \to +\infty]{} 0.$$
Then,
$$\E[f(X^n_{t_1}, \ldots , X^n_{t_k})g(\tau^{*,n})] - \E[f(X_{t_1}, \ldots , X_{t_k})g(\tau^*)] \xrightarrow[n \to +\infty]{} 0.$$
Let us now consider $\varphi : \R^{k+1} \to \R$ continuous and bounded. \\
Let us fix $\varepsilon > 0$. \\
$((X^n_{t_1}, \ldots , X^n_{t_k},\tau^{*,n}))_n$ is tight. We can find a compact set $K_\varepsilon$ such that 
\begin{equation}
\label{beta1}
\P[(X^n_{t_1}, \ldots , X^n_{t_k},\tau^{*,n}) \notin K_\varepsilon] \leqslant \varepsilon.
\end{equation}
We write $\varphi=\varphi 1_{K_\varepsilon} + \varphi 1_{K_\varepsilon^c}$. \\
$\varphi 1_{K_\varepsilon}$ is a continuous function on the compact set $K_\varepsilon$. Using Weierstrass' Theorem, we can find a
polynomial function $P$ such that
\begin{equation}
\label{beta2}
\|\varphi 1_{K_\varepsilon} - P 1_{K_\varepsilon} \|_{\infty} \leqslant \varepsilon.
\end{equation}
Using the previous result and the linearity of expectation, we have 
\begin{equation}
\label{beta3}
\E[(P 1_{K_\varepsilon})(X^n_{t_1}, \ldots , X^n_{t_k},\tau^{*,n})] - \E[(P 1_{K_\varepsilon})(X_{t_1}, \ldots , X_{t_k},\tau^*)] 
        \xrightarrow[n \to +\infty]{} 0.
\end{equation}
Finally, 
\begin{eqnarray*}
\lefteqn{ |\E[\varphi(X^n_{t_1}, \ldots , X^n_{t_k},\tau^{*,n})]-\E[\varphi(X_{t_1}, \ldots , X_{t_k},\tau^*)]|} \\
& = & |\E[(\varphi(X^n_{t_1}, \ldots , X^n_{t_k},\tau^{*,n})]-\E[\varphi(X_{t_1}, \ldots , X_{t_k},\tau^*))
        1_{K_\varepsilon}(X^n_{t_1}, \ldots , X^n_{t_k},\tau^{*,n})]| \\
& &  + |\E[(\varphi(X^n_{t_1}, \ldots , X^n_{t_k},\tau^{*,n})]-\E[\varphi(X_{t_1}, \ldots , X_{t_k},\tau^*))
        1_{K_\varepsilon^c}(X^n_{t_1}, \ldots , X^n_{t_k},\tau^{*,n})]| \\
& \leqslant & 2 \|\varphi 1_{K_\varepsilon} - P 1_{K_\varepsilon} \|_{\infty} \\
& & \quad + \E[(P 1_{K_\varepsilon})(X^n_{t_1}, \ldots , X^n_{t_k},\tau^{*,n})] - \E[(P 1_{K_\varepsilon})(X_{t_1}, \ldots , X_{t_k},\tau^*)] \\
& & \quad \quad + \| \varphi \|_\infty \P[(X^n_{t_1}, \ldots , X^n_{t_k},\tau^{*,n}) \notin K_\varepsilon] \\
& \leqslant & \E[(P 1_{K_\varepsilon})(X^n_{t_1}, \ldots , X^n_{t_k},\tau^{*,n})] - \E[(P 1_{K_\varepsilon})(X_{t_1}, \ldots , X_{t_k},\tau^*)]\\
& & \quad + (2 + \| \varphi \|_\infty) \varepsilon \text{~~using (\ref{beta1}) and (\ref{beta2}).}
\end{eqnarray*}
Taking the limit for $n$, using (\ref{beta3}), we obtain :
$$\limsup_n |\E[\varphi(X^n_{t_1}, \ldots , X^n_{t_k},\tau^{*,n})]-\E[\varphi(X_{t_1}, \ldots , X_{t_k},\tau^*)]| 
        \leqslant (2 + \| \varphi \|_\infty) \varepsilon. $$
This is true for every $\varepsilon >0$, so we have 
$$\E[\varphi(X^n_{t_1}, \ldots , X^n_{t_k},\tau^{*,n})] \xrightarrow[n \to +\infty]{} \E[\varphi(X_{t_1}, \ldots , X_{t_k},\tau^*)].$$
Then, $(X^n_{t_1}, \ldots , X^n_{t_k},\tau^{*,n}) \cvl (X_{t_1}, \ldots , X_{t_k}, \tau^*)$. \\
\indent The tightness of the sequence $((X^n, \tau^{*,n}))_n$ and the finite-dimensional convergence on a dense set to $(X, \tau^*)$ implies 
$(X^n, \tau^{*,n}) \cvl (X,\tau^*)$.
\findem

\subsection{Application to the proof of the inequality $\limsup \Gamma_n(L) \leqslant \Gamma(L)$}

We can now prove a result about the convergence of optimal values.

\begin{theo}
\label{thlimsup1}
Let us consider a c\`adl\`ag process X continuous in probability, its natural right continuous filtration $\F$, a sequence $(X^n)_n$ of c\`adl\`ag 
processes and their natural filtrations $(\F^n)_n$. We suppose that $X^n \cvp X$ and $\forall n$, $\F^n \subset \F$. \\
Then $\limsup \Gamma_n(L) \leqslant \Gamma(L)$.
\end{theo}

\demo
There exists a subsequence $(\Gamma_{\varphi(n)}(L))_n$ converging to $\limsup \Gamma_n(L)$. \\
Let us fix $\varepsilon > 0$. We can find a sequence $(\tau^{\varphi(n)})_n$ of $(\mathcal{T}^{\varphi(n)}_L)_n$ such that 
$$\forall n, \E[\gamma(\tau^{\varphi(n)}, X^{\varphi(n)}_{\tau^{\varphi(n)}})] \geqslant \Gamma_{\varphi(n)}(L) - \varepsilon.$$
We consider the sequence $(\tau^{*,n})_n$ of randomized stopping times associated to \mbox{$(\tau^{n})_n$ :} for every $n$, 
$\tau^{*,n}(\omega,t)=\tau^{n}(\omega)$, $\forall \omega$, $\forall t$. \\
$\F^n \subset \F$ and $(\tau^{\varphi(n)})$ is a sequence of $(\F^{\varphi(n)})_n$-stopping times bounded by $L$, so using Proposition 
\ref{sssuitetar}, there exists a randomized $\F$ stopping time $\tau^*$ and a 
subsequence $(\tau^{\varphi \circ \psi(n)})$ such that $\tau^{*,\varphi \circ \psi(n)} \xrightarrow{BC} \tau^*$. \\
$X^{\varphi \circ \psi(n)} \cvp X$ and $\tau^{*,\varphi \circ \psi(n)} \xrightarrow{BC} \tau^*$, so using Proposition \ref{cvcouple}, 
$$(X^{\varphi \circ \psi(n)}, \tau^{*,\varphi \circ \psi(n)}) \cvl (X,\tau^*).$$
Then, using Proposition \ref{cvlXntar}, we have :
$$(\tau^{*,\varphi \circ \psi (n)},X^{\varphi \circ \psi (n)}_{\tau^{*,\varphi \circ \psi (n)}}) \cvl (\tau^*,X_{\tau^*}).$$
Since $\gamma$ is continuous and bounded, we have :
$$\E[\gamma(\tau^{*,\varphi \circ \psi (n)}, X^{\varphi \circ \psi(n)}_{\tau^{*,\varphi \circ \psi(n)}})] \to \E[\gamma(\tau^*,X_{\tau^*})].$$
But,  
$\E[\gamma(\tau^{*,\varphi \circ \psi(n)}, X^{\varphi \circ \psi(n)}_{\tau^{*,\varphi \circ \psi(n)}})] 
	= \E[\gamma(\tau^{\varphi \circ \psi(n)},X^{\varphi \circ \psi(n)}_{\tau^{\varphi \circ \psi(n)}})]$ 
by definition of $(\tau^{*,n})$ and by choice of $\varphi$, 
$\E[\gamma(\tau^{\varphi \circ \psi(n)},X^{\varphi \circ \psi(n)}_{\tau^{\varphi \circ \psi(n)}})] 
	\geqslant \Gamma_{\varphi \circ \psi(n)}(L) - \varepsilon$.
So, 
$$\E[\gamma(\tau^*,X_{\tau^*})] \geqslant \limsup \Gamma_{\varphi \circ \psi(n)}(L) - \varepsilon.$$
By selection of $\varphi$, $\limsup \Gamma_{\varphi \circ \psi(n)}(L)=\limsup \Gamma_{n}(L).$ Then, 
$$\E[\gamma(\tau^*,X_{\tau^*})] \geqslant \limsup \Gamma_{n}(L) - \varepsilon.$$
This is true for every $\varepsilon >0$, so
$$\E[\gamma(\tau^*,X_{\tau^*})] \geqslant \limsup \Gamma_n(L).$$
But, by definition, $\E[\gamma(\tau^*,X_{\tau^*})] \leqslant \Gamma^*(L)$ because 
$\tau^*$ is a randomized stopping time. 
As $\Gamma^*(L)=\Gamma(L)$ by Lemma \ref{GammaGamma*}, we deduce $\Gamma(L) \geqslant \limsup \Gamma_n(L)$.
\findem

\begin{rem} 
In the previous Theorem, the most important argument is that we know things about the nature of the limit of the subsequence of stopping
times thanks to Proposition \ref{sssuitetar}. If we remove the inclusion of the filtrations $\F^n \subset \F, \forall n$, the limit of the 
subsequence is no longer a randomized $\F$ stopping time. In this case, we can't compare $\E[\gamma(\tau^*,X_{\tau^*})]$ and $\Gamma^*(L)$.
\end{rem}


\section{Proof of the inequality $\limsup \Gamma_n(L) \leqslant \Gamma(L)$ when $\F^n \xrightarrow{w} \F$}
\label{limsup_cvfiltrations}

\begin{theo}
\label{thlimsup2}
Let us consider a sequence of c\`adl\`ag processes $(X^n)_n$, their na\-tu\-ral filtrations $(\F^n)_n$, a c\`adl\`ag process continuous in 
probability $X$ and its right continuous natural filtration $\F$. We suppose $X^n \cvp X$ and $\F^n \xrightarrow{w} \F$. \\
Then $\limsup \Gamma_n(L) \leqslant \Gamma(L)$.
\end{theo}

\demo
We argue more or less as Aldous in the second part of the proof of \citep[Theorem 17.2]{preprintAldous}. \\
We can find a subsequence $(\Gamma_{\varphi(n)}(L))_n$ converging to $\limsup \Gamma_n(L)$. \\
Let us take $\varepsilon > 0$. There exists a sequence $(\tau^{\varphi(n)})_n$ of $(\mathcal{T}^{\varphi(n)}_L)_n$ such that 
$$\forall n, \E[\gamma(\tau^{\varphi(n)}, X^{\varphi(n)}_{\tau^{\varphi(n)}})] \geqslant \Gamma_{\varphi(n)}(L) - \varepsilon.$$
Let us consider the sequence $(\tau^{*,n})_n$ of associated randomized $(\F^n)$ stopping times like in \ref{tar}. Taking the filtration 
$\mathcal{H}=(\bigvee_n \F^n) \vee \F$, $(\tau^{*,n})$ is a bounded sequence of randomized $\mathcal{H}$ stopping times. Then,
using \citep[Theorem 1.5]{BC}, we can find an increasing map $\varphi$ and a randomized $\mathcal{H}$
stopping time $\tau^*$ ($\tau^*$ is not a priori a randomized $\F$ stopping time) such that 
$$\tau^{*,\varphi(n)} \cvBC \tau^*.$$
Using Proposition \ref{cvcouple}, we obtain $(X^{\varphi(n)}, \tau^{*,\varphi(n)}) \cvl (X,\tau^*)$.
Then, with Proposition \ref{cvlXntar}, we have
$(\tau^{\varphi(n)},X^{\varphi(n)}_{\tau^{\varphi(n)}}) \cvl (\tau^*,X_{\tau^*}).$ So,
$$\E[\gamma(\tau^{\varphi(n)},X^{\varphi(n)}_{\tau^{\varphi(n)}})] 
        \xrightarrow[n \to +\infty]{} \E[\gamma(\tau^*,X_{\tau^*})].$$
On the other hand, $\E[\gamma(\tau^{\varphi(n)},X^{\varphi(n)}_{\tau^{\varphi(n)}})]
        \geqslant \Gamma_{\varphi(n)}(L) - \varepsilon$. So, when $n$ tends to infinity, it results :
$$\E[\gamma(\tau^*,X_{\tau^*})] \geqslant \limsup \Gamma_n(L) - \varepsilon.$$
This is true for every $\varepsilon >0$, hence we have 
\begin{equation}
\label{inequgamma}
\E[\gamma(\tau^*,X_{\tau^*})] \geqslant \limsup \Gamma_n(L).
\end{equation}
It remains to compare $\E[\gamma(\tau^*,X_{\tau^*})]$ and $\Gamma(L)$. \\
Let us consider the smaller right continuous filtration $\G$ such that $X$ is $\G$ adapted and $\tau^*$ is a randomized $\G$ stopping time. 
It is clear that $\F \subset \G$. 
For every $t$, we have 
$$\G_t \times \mathcal{B}=\cap_{s > t} \sigma(A \times [0,1] ,\{\tau^* \leqslant u\}, A \in \F_s, u \leqslant s).$$ 

We consider the set $\tilde{\mathcal{T}}_L$ of randomized $\G$ stopping times bounded by $L$ and we define 
$\tilde{\Gamma}(L)=\underset{\tilde{\tau} \in \tilde{\mathcal{T}}_L}{\sup}\E[\gamma(\tilde{\tau},X_{\tilde{\tau}})]$. \\

By definition of $\G$, $\tau^* \in \tilde{\mathcal{T}}_L$ so $\E[\gamma(\tau^*,X_{\tau^*})] \leqslant \tilde{\Gamma}(L)$. \\ 

We are going to end the proof using the following Lemma, that is an adaptation of \citep[Proposition 3.5]{cvreduites} to our enlargement 
of filtration :
\begin{lemme}
\label{indepcond}
If $\G_t \times \mathcal{B}$ and $\F_T \times \mathcal{B}$ are conditionally independent given $\F_t \times \mathcal{B}$ 
for every $t \in [0,T]$, then $\tilde{\Gamma}(L)=\Gamma^*(L)$.
\end{lemme}

\demo
The proof is the same as the proof of \citep[Proposition 3.5]{cvreduites} with 
$(\F_t \times \B)_{t \in [0,T]}$ and $(\G_t \times \B)_{t \in [0,T]}$ instead of $\F^Y$ and $\F$ and with the process $X^*$ such that 
for every $\omega$, for every $v \in [0,1]$, for every $t \in [0,T]$, 
$X^*_t(\omega,v)=X_t(\omega)$ instead of the process $Y$.
\findem

According to \citep[Theorem 3]{BY}, the condition of conditional independence required in Lemma
\ref{indepcond} is equivalent to the following assumption : 
\begin{equation}
\label{indcond}
\forall t \in [0,T], \forall Z \in L^1(\F_T \times \mathcal{B}), \E[Z|\F_t \times \mathcal{B}]=\E[Z|\G_t \times \mathcal{B}].
\end{equation}

We will show that the assumptions of Theorem \ref{thlimsup2} imply those of Lemma \ref{indepcond}, therefore proving inequality 
(\ref{indcond}). \\
Note that in \citep{preprintAldous} and in \citep{cvreduites}, they need extended convergence to prove (\ref{indcond}). \\

\noindent Without loss of generality, we suppose from now that $\tau^n \cvBC \tau$ instead of $\tau^{\varphi(n)} \cvBC \tau$. \\
We also denote by "continuity points" of a process the points where the process is continuous in probability, ie $t$ such that 
$\P[\Delta X_t \not= 0]=0$.\\

\noindent - As $\F \subset \G$, for every $t$, $\forall Z \in L^1(\F_T \times \mathcal{B})$, 
$\E_{\P \otimes \mu}[Z|\F_t \times \mathcal{B}]$ is $\G_t \times \mathcal{B}$-measurable. \\
- Let us show $\forall t \in [0,T], \forall Z \in L^1(\F_T \times \mathcal{B}), \forall C \in \G_t \times \mathcal{B},$ 
$$\E_{\P \otimes \mu}[\E_{\P \otimes \mu}[Z|\F_t \times \mathcal{B}]1_C]=\E_{\P \otimes \mu}[Z1_C].$$

Let us fix $t \in [0,T]$ and $\varepsilon > 0$. \\
Let us take $Z \in L^1(\F_T \times \mathcal{B})$. 
By definition of $\G_t \times \mathcal{B}$, it suffices to prove that for every $A \in \F_t$, for every $s \leqslant t$ and for every 
$B \in \mathcal{B}$, 
\begin{eqnarray}
\label{eqZ}
\lefteqn{\int\int_{\Omega \times [0,1]} Z(\omega,v)1_A(\omega)1_{\{\tau^*(\omega,v) \leqslant s\}}1_B(v) d\P(\omega)dv }\\
& = &
\int\int_{\Omega \times [0,1]} \E_{\P \otimes \mu}[Z|\F_t \times \mathcal{B}](\omega,v)1_A(\omega)1_{\{\tau^*(\omega,v) \leqslant s\}}
                1_B(v) d\P(\omega)dv. \nonumber
\end{eqnarray}

We first prove that (\ref{eqZ}) holds for $Z=1_{A_1 \times A_2}$, $A_1 \in \F_T$, $A_2 \in \mathcal{B}$. \\
We can find $l \in \N$, some continuity points of $X$ $s_1 < \ldots < s_l$ and a continuous bounded function $f$ such that  
\begin{equation}
\label{alpha1}
\E_{\P}[|1_{A_1}-f(X_{s_1}, \ldots, X_{s_l})|] \leqslant \varepsilon.
\end{equation}
Then
$$\int\int |1_{A_1 \times A_2}(\omega,v)-f(X_{s_1}(\omega), \ldots, X_{s_l}(\omega))1_{A_2}(v)|d\P(\omega)dv \leqslant \varepsilon.$$

Let us fix $A \in \F_t$. 
We can find $k \in \N$, $t_1 < \ldots < t_k \leqslant t$ where $t_i$ are continuity points of $X$ and $H:\R^k \to \R$ bounded continuous such
that 
\begin{equation}
\label{alpha2}
\E_{\P}[|1_A-H(X_{t_1}, \ldots, X_{t_k})|] \leqslant \varepsilon.
\end{equation}
Let $u \geqslant t$ be a continuity point of $\E[f(X_{s_1}, \ldots, X_{s_l})|\F_.]$ and of $\tau^*$. \\
Fix $s \leqslant t$. 
We can find $G$ bounded continuous such that 
\begin{equation}
\label{alpha3}
\E_{\P \otimes \mu}[|1_{\{\tau^* \leqslant s\}}-G(\tau^* \wedge u)|] \leqslant \varepsilon.
\end{equation}
$B \in \mathcal{B}$ and the set of continuous functions is dense into $L^1(\mu)$, so there exists $g : \R \to \R$ bounded continuous such
that 
\begin{equation}
\label{alpha45}
\int |1_{B}(v)-g(v)|dv \leqslant \varepsilon.
\end{equation}

\noindent We are going to show that  
\begin{eqnarray*}
&& \int \int \E_{\P \otimes \mu}[f(X_{s_1}, \ldots, X_{s_l})1_{A_2}|\F_u \otimes \mathcal{B}](\omega,v)
                H(X_{t_1}(\omega), \ldots, X_{t_k}(\omega))\\
&& \quad \quad \quad G(\tau^*(\omega,v) \wedge u)g(v)d(\P \otimes \mu)(\omega,v)\\
& = & \int\int f(X_{s_1}(\omega), \ldots, X_{s_l}(\omega))1_{A_2}(v)H(X_{t_1}(\omega), \ldots, X_{t_k}(\omega))\\
&& \quad \quad \quad G(\tau^*(\omega,v) \wedge u)g(v)d(\P \otimes \mu)(\omega,v).
\end{eqnarray*}
$X^n \cvp X$, $s_i$ are continuity points of $X$ and $f$ is a bounded continuous function, then
\begin{equation}
\label{mer0} 
f(X^n_{s_1}, \ldots, X^n_{s_l}) \cvL1 f(X_{s_1}, \ldots, X_{s_l}).
\end{equation}
Moreover, $\F^n \xrightarrow{w} \F$ so using \citep[Remark 2]{cvfiltration}, 
$$\E_{\P}[f(X^n_{s_1}, \ldots, X^n_{s_l})|\F^n] \cvp \E_{\P}[f(X_{s_1}, \ldots, X_{s_l})|\F].$$ 
Since $u$ is a continuity point of $\E_{\P}[f(X_{s_1}, \ldots, X_{s_l})|\F_.]$, we have 
$$\E_{\P}[f(X^n_{s_1}, \ldots, X^n_{s_l})|\F^n_u] \cvp \E_{\P}[f(X_{s_1}, \ldots, X_{s_l})|\F_u].$$
Since $f$ is bounded, convergence is in $L^1$ :
\begin{equation}
\label{mer1}
\E_{\P}[f(X^n_{s_1}, \ldots, X^n_{s_l})|\F^n_u] \cvL1 \E_{\P}[f(X_{s_1}, \ldots, X_{s_l})|\F_u].
\end{equation}
Let us consider the maps $\tilde{H}$, $\tilde{G}$ and $\tilde{g}$ from $\R^{k+l+2}$ to $\R$ defined as follows : \\
\begin{displaymath}
\begin{array}{rcl}
\tilde{H}(x_1, \ldots , x_l, y_1, \ldots , y_k, z, v) & = & H(y_1, \ldots , y_k), \\
\tilde{G}(x_1, \ldots , x_l, y_1, \ldots , y_k, z, v) & = & G(z),\\
\tilde{g}(x_1, \ldots , x_l, y_1, \ldots , y_k, z, v) & = & g(v).
\end{array}
\end{displaymath}
Fix $\varepsilon'>0$. \\
The set of continuous maps is dense into \mbox{$L^1(\P_{(X_{s_1}, \ldots, X_{s_l})} \otimes \mu)$,} 
hence we can find $h : \R^{l+1} \to \R$ such that
\begin{equation}
\label{mer25}
\int\int |h(x_1, \ldots , x_l, v) - f(x_1, \ldots , x_l)1_{A_2}(v)|d(\P_{(X_{s_1}, \ldots, X_{s_l})} \otimes \mu)(\omega,v) 
\leqslant \varepsilon'.
\end{equation}
$X^n \cvp X$, $s_i$ are continuity points of $X$ and $h$ is a bounded continuous function, so
\begin{equation}
\label{mer27} 
\int\int|h(X^n_{s_1}(\omega), \ldots, X^n_{s_l}(\omega), v) - h(X_{s_1}(\omega), \ldots, X_{s_l}(\omega), v)|d(\P \otimes \mu)(\omega,v) 
\xrightarrow[n \to +\infty]{} 0.
\end{equation}
Then we consider :
$$\tilde{h}(x_1, \ldots , x_l, y_1, \ldots , y_k, z, v) = h(x_1, \ldots , x_l, v). $$
$\tilde{h}\tilde{H}\tilde{G}\tilde{g}$ is continuous as product of continuous maps. \\
Moreover, $(X^n,\tau^{*,n}) \cvl (X,\tau^*)$ and $u$ is a continuity point of $\tau^*$, so that $(X^n,\tau^{*,n} \wedge u) \cvl (X,\tau^* \wedge u)$. \\
Let $U : \Omega \times [0,1] \to [0,1]$ be the random variable such that $\forall \omega, \forall v, U(\omega,v)=v$. As in the proof of 
Proposition \ref{cvcouple}, we have :
$$(X^n,\tau^{*,n} \wedge u, U) \cvl (X,\tau^* \wedge u, U).$$
As $s_1, \ldots, s_l, t_1, \ldots, t_k$ are continuity points of $X$, we have 
$$(X^n_{s_1}, \ldots , X^n_{s_l}, X^n_{t_1}, \ldots , X^n_{t_k}, \tau^{*,n} \wedge u, U)
\cvl (X_{s_1}, \ldots , X_{s_l}, X_{t_1}, \ldots , X_{t_k}, \tau^* \wedge u, U).$$
Hence, 
\begin{eqnarray}
\label{mer26}
\lefteqn{\E_{\P \otimes \mu}[
        (\tilde{h}\tilde{H}\tilde{G}\tilde{g})(X^n_{s_1}, \ldots , X^n_{s_l}, X^n_{t_1}, \ldots , X^n_{t_k}, \tau^{*,n} \wedge u, U)] } \\
&\xrightarrow[n \to +\infty]{} & 
\E_{\P \otimes \mu}[(\tilde{h}\tilde{H}\tilde{G}\tilde{g})(X_{s_1}, \ldots , X_{s_l}, X_{t_1}, \ldots , X_{t_k},\tau^* \wedge u, U)]. 
\nonumber
\end{eqnarray}
By definition of functions $\tilde{h}$, $\tilde{H}$, $\tilde{G}$ and $\tilde{g}$, we have :
\begin{eqnarray}
\label{cv1}
&&\int\int h(X^{n}_{s_1}(\omega), \ldots, X^{n}_{s_l}(\omega),v)H(X^{n}_{t_1}(\omega), \ldots , X^{n}_{t_k}(\omega))\nonumber\\
&& \quad \quad \quad \quad G(\tau^{*,n}(\omega,v) \wedge u)g(v)d(\P \otimes \mu)(\omega,v)   \nonumber\\
& \xrightarrow[n \to +\infty]{} & 
\int\int h(X_{s_1}(\omega), \ldots, X_{s_l}(\omega), v)H(X_{t_1}(\omega), \ldots , X_{t_k}(\omega))\nonumber\\
&& \quad \quad \quad \quad G(\tau^*(\omega,v) \wedge u)g(v)d(\P \otimes \mu)(\omega,v).
\end{eqnarray}
Then using triangular inequalities and inequations (\ref{mer0}), (\ref{mer27}), (\ref{cv1}) and (\ref{mer25}), we get
\begin{eqnarray*}
& & \limsup_n \bigg|\int\int f(X^{n}_{s_1}(\omega), \ldots, X^{n}_{s_l}(\omega))1_{A_2}(v)H(X^{n}_{t_1}(\omega), \ldots ,X^{n}_{t_k}(\omega)) \\
& & \quad \quad         G(\tau^{*,n}(\omega,v) \wedge u)g(v)d(\P \otimes \mu)(\omega,v) \\ 
& & - \int\int f(X_{s_1}(\omega), \ldots, X_{s_l}(\omega))1_{A_2}(v)H(X_{t_1}(\omega), \ldots , X_{t_k}(\omega))\\
& & \quad \quad G(\tau^*(\omega,v) \wedge u)g(v)d(\P \otimes \mu)(\omega,v)\bigg|\\
& \leqslant &  2\|H\|_\infty \|G\|_\infty \|g\|_\infty \varepsilon'.
\end{eqnarray*}
This is true for every $\varepsilon' >0$, so :
\begin{eqnarray}
\label{mer28}
&&\int\int f(X^{n}_{s_1}(\omega), \ldots, X^{n}_{s_l}(\omega))1_{A_2}(v)H(X^{n}_{t_1}(\omega), \ldots ,X^{n}_{t_k}(\omega)) \nonumber\\
&&      \quad \quad \quad \quad         G(\tau^{*,n}(\omega,v) \wedge u)g(v)d(\P \otimes \mu)(\omega,v)  \\
&\xrightarrow[n \to +\infty]{} & 
\int\int f(X_{s_1}(\omega), \ldots, X_{s_l}(\omega))1_{A_2}(v)H(X_{t_1}(\omega), \ldots , X_{t_k}(\omega))\nonumber\\
&&      \quad \quad \quad \quad G(\tau^*(\omega,v) \wedge u)g(v)d(\P \otimes \mu)(\omega,v). 
\nonumber
\end{eqnarray}

\medskip

On the other hand, $\E[f(X_{s_1}, \ldots, X_{s_l})1_{A_2}|\F_u \times \mathcal{B}]=\E[f(X_{s_1}, \ldots, X_{s_l})|\F_u]1_{A_2}$. \\
$\E[f(X_{s_1}, \ldots, X_{s_l})|\F_u]$ is $\F_u$-measurable. \\
Let us fix $\varepsilon' > 0.$ We can find $j \in \N$ and $v_1 < \ldots < v_j \leqslant u$ some continuity points of $X$ and $F:\R^j \to \R$
bounded continuous such that :
\begin{equation}
\label{mer2}
\E_{\P}[|\E[f(X_{s_1}, \ldots, X_{s_l})|\F_u] - F(X_{v_1}, \ldots, X_{v_j})|] \leqslant \varepsilon'.
\end{equation}
$X^n \cvp X$, $F$ is bounded continuous and $v_i$ are continuity points of $X$ then 
\begin{equation}
\label{mer4}
F(X^n_{v_1}, \ldots, X^n_{v_j}) \cvL1 F(X_{v_1}, \ldots, X_{v_j}).
\end{equation} 
As previously, we have :
\begin{eqnarray}
\label{cv2.1}
&& \int\int F(X^n_{v_1}(\omega), \ldots, X^n_{v_j}(\omega))1_{A_2}(v)H(X^{n}_{t_1}(\omega), \ldots , X^{n}_{t_k}(\omega))\nonumber\\
&& \quad \quad \quad \quad G(\tau^{*,n}(\omega,v) \wedge u)g(v)d(\P \otimes \mu)(\omega,v)   \nonumber\\
& \xrightarrow[n \to +\infty]{} & 
\int\int F(X_{v_1}(\omega), \ldots, X_{v_j}(\omega))1_{A_2}(v)H(X_{t_1}(\omega), \ldots , X_{t_k}(\omega))\nonumber\\
&& \quad \quad \quad \quad G(\tau^*(\omega, v) \wedge u)g(v)d(\P \otimes \mu)(\omega,v).
\end{eqnarray}
Then using triangular inequalities and inequations (\ref{mer0}), (\ref{mer1}), (\ref{mer2}), (\ref{mer4}) and (\ref{cv2.1}), we have :
\begin{displaymath}
\begin{array}{l}
 \limsup_n \bigg| 
\int\int \E[f(X^n_{s_1}, \ldots, X^n_{s_l})1_{A_2}|\F^n_u \times \mathcal{B}](\omega, v)
        H(X^{n}_{t_1}(\omega), \ldots , X^{n}_{t_k}(\omega))\\
\quad \quad \quad \quad G(\tau^{*,n} \wedge u)g(v)d(\P \otimes \mu)(\omega, v) \\
 \quad - \int\int \E[f(X_{s_1}, \ldots, X_{s_l})1_{A_2}|\F_u \times \mathcal{B}](\omega, v)
        H(X_{t_1}(\omega), \ldots , X_{t_k}(\omega))\\
\quad \quad \quad \quad	G(\tau^*(\omega, v) \wedge u)g(v)d(\P \otimes \mu)(\omega, v) \bigg|  \\
 \leqslant  2\|H\|_{\infty}\|G\|_{\infty}\|g\|_{\infty} \varepsilon'.
\end{array}
\end{displaymath}
This is true for every $\varepsilon' > 0$, so :
\begin{eqnarray}
\label{cv2}
&&\int\int \E[f(X^n_{s_1}, \ldots, X^n_{s_l})1_{A_2}|\F^n_u \times \mathcal{B}](\omega, v)
        H(X^{n}_{t_1}(\omega), \ldots , X^{n}_{t_k}(\omega))\nonumber\\
&& \quad \quad 	\quad \quad	G(\tau^{*,n} \wedge u)g(v)d(\P \otimes \mu)(\omega, v)\nonumber \\
& \xrightarrow[n \to +\infty]{} & 
\int\int \E[f(X_{s_1}, \ldots, X_{s_l})1_{A_2}|\F_u \times \mathcal{B}](\omega, v)
        H(X_{t_1}(\omega), \ldots , X_{t_k}(\omega))\nonumber \\
&& \quad \quad  \quad \quad	G(\tau^*(\omega, v) \wedge u)g(v)d(\P \otimes \mu)(\omega, v). 
\end{eqnarray}

\medskip

But, $H(X^{n}_{t_1}, \ldots , X^{n}_{t_k})$ is $\F^{n}_u \times \mathcal{B}$-measurable and 
\mbox{$G(\tau^{n} \wedge u)$} and $g(U)$ are also $\F^{n}_u \times \mathcal{B}$-measurable, by continuity of $G$ and $g$. Then,
\begin{eqnarray*}
\lefteqn{\E[\E[f(X^{n}_{s_1}, \ldots, X^{n}_{s_l})1_{A_2}|\F^{n}_u \times \mathcal{B}]H(X^{n}_{t_1}, \ldots , X^{n}_{t_k})G(\tau^{n} \wedge u)g(U)] }\\
& = & \E[\E[f(X^{n}_{s_1}, \ldots, X^{n}_{s_l})1_{A_2}H(X^{n}_{t_1}, \ldots , X^{n}_{t_k})G(\tau^{n} \wedge u)g(U)|\F^{n}_u \times \mathcal{B}]] \\
& = & \E[f(X^{n}_{s_1}, \ldots, X^{n}_{s_l})1_{A_2}H(X^{n}_{t_1}, \ldots , X^{n}_{t_k})G(\tau^{n} \wedge u)g(U)]
\end{eqnarray*}
Using unicity of the limit and convergences (\ref{cv1}) and (\ref{cv2}), we obtain :
\begin{eqnarray}
\label{eqfHG}
&& \int\int \E[f(X_{s_1}, \ldots, X_{s_l})1_{A_2}|\F_u \times \mathcal{B}](\omega, v)
        H(X_{t_1}(\omega), \ldots , X_{t_k}(\omega))\nonumber \\
&& \quad \quad  \quad \quad	G(\tau^*(\omega, v) \wedge u)g(v)d(\P \otimes \mu)(\omega, v) \nonumber \\
& = & \int\int f(X_{s_1}(\omega), \ldots, X_{s_l}(\omega))1_{A_2}(v)
        H(X_{t_1}(\omega), \ldots , X_{t_k}(\omega))\nonumber \\ 
&& \quad \quad  \quad \quad G(\tau^*(\omega, v) \wedge u)g(v)d(\P \otimes \mu)(\omega, v).
\end{eqnarray}
Then, 
\begin{displaymath}
\begin{array}{l}
 \Big| \int\int \E[f(X_{s_1}, \ldots, X_{s_l})1_{A_2}|\F_u \times \mathcal{B}](\omega, v)1_A(\omega)
        1_{\{\tau^*(\omega, v) \leqslant s\}}1_{B}(v)d(\P \otimes \mu)(\omega, v)\\
 \quad - \int\int f(X_{s_1}(\omega), \ldots, X_{s_l}(\omega))1_{A_2}(v)1_A(\omega)
        1_{\{\tau^*(\omega, v) \leqslant s\}}1_{B}(v)d(\P \otimes \mu)(\omega, v) \Big| \\
\leqslant 
 2\|f\|_\infty (1 + \|H\|_\infty + \|G\|_\infty) \varepsilon 
 \text{~~using (\ref{alpha2}), (\ref{alpha3}), (\ref{alpha45}) and (\ref{eqfHG})}.
\end{array}
\end{displaymath}
Let $u$ tend to $t$ by upper values. $\E[f(X_{s_1}, \ldots, X_{s_l})|\F_.]$ is a c\`adl\`ag process, so we have :
\begin{eqnarray}
\label{alpha4}
& \Big| \int\int \E[f(X_{s_1}, \ldots, X_{s_l})1_{A_2}|\F_t \times \mathcal{B}](\omega, v)1_A(\omega)
        1_{\{\tau^*(\omega, v) \leqslant s\}}1_B(v)d(\P \otimes \mu)(\omega, v) & \nonumber \\
&\quad - \int\int f(X_{s_1}(\omega), \ldots, X_{s_l}(\omega))1_{A_2}(v)1_A(\omega) 
        1_{\{\tau^*(\omega, v) \leqslant s\}}1_B(v)d(\P \otimes \mu)(\omega, v) \Big| &
\nonumber \\
& \leqslant  2\|f\|_\infty (1 + \|H\|_\infty + \|G\|_\infty) \varepsilon. &
\end{eqnarray}
Then, 
\begin{eqnarray*}
&& \Big| \int\int \E[Z|\F_t \times \mathcal{B}](\omega, v)1_A(\omega)1_{\{\tau^*(\omega, v) \leqslant s\}}d(\P \otimes \mu)(\omega, v)\\
&& \quad - \int\int Z(\omega, v)1_A(\omega)1_{\{\tau^*(\omega, v) \leqslant s\}}1_B(v)d(\P \otimes \mu)(\omega, v) \Big| \\
& \leqslant & 
2\|f\|_\infty (1 + \|H\|_\infty + \|G\|_\infty) \varepsilon + 2\varepsilon \text{~~using (\ref{alpha1}) and (\ref{alpha4})}.
\end{eqnarray*}
This is true for every $\varepsilon > 0$, so we have the equality (\ref{eqZ}) :
\begin{eqnarray*}
\lefteqn{
\int\int \E[Z|\F_t \times \mathcal{B}](\omega, v)1_A(\omega)1_{\{\tau^*(\omega, v) \leqslant s\}}1_B(v)d(\P \otimes \mu)(\omega,v)}\\
& = & \int\int Z(\omega, v)1_A(\omega)1_{\{\tau^*(\omega, v) \leqslant s\}}1_B(v)d(\P \otimes \mu)(\omega, v), 
\end{eqnarray*}
for every $t \in [0,T]$, for every $Z =1_{A_1 \times A_2}$, $A_1 \in \F_T$, $A_2 \in \mathcal{B}$, for every $A \in \F_t$, 
for every $s \leqslant t$, for every $B \in \mathcal{B}$. \\

If $Z=1_E$ with $E \in \F_T \times \mathcal{B}$, (\ref{eqZ}) holds using the preceding results and an argument of monotone class. \\

Then, if $Z$ is a function of the form $\sum a_i 1_{A_i}$ with $a_i \in \R$ and $A_i \in \F_T \times \mathcal{B}$, (\ref{eqZ}) holds by
linearity. \\

If $Z$ is $\F_T \times \mathcal{B}$-measurable, we use density in $L^1$ norm of the functions of the form $\sum a_i 1_{A_i}$ to obtain 
(\ref{eqZ}). \\

Hence, for every $t \in [0,T]$, for every $Z \in L^1(\F_T \times \mathcal{B})$, 
for every $C \in \G_t \times \mathcal{B}$ (by definition of $\G_t \times \mathcal{B}$),
$$\E_{\P \otimes \mu}[\E_{\P \otimes \mu}[Z|\F_t \times \mathcal{B}]1_C]=\E_{\P \otimes \mu}[Z1_C].$$

\noindent The assumption of Lemma \ref{indepcond} if filled, so 
$$\E[\gamma(\tau,Y_\tau)] \leqslant \tilde{\Gamma}(L)=\Gamma^*(L).$$
Using inequality (\ref{inequgamma}), we finally have
$$\limsup \Gamma_n(L) \leqslant \Gamma^*(L).$$
But using Lemma \ref{GammaGamma*}, $\Gamma^*(L)=\Gamma(L)$.
Theorem \ref{thlimsup2} is proved.
\findem

To sum up, under the hypothesis of Theorem \ref{cvGamma}, we have proved the inequality $\Gamma(L) \leqslant \liminf \Gamma_n(L)$ in Theorem 
\ref{thliminf} and Remark \ref{rqF}. Then, we have shown that $\Gamma(L) \geqslant \liminf \Gamma_n(L)$ when we have the inclusion of 
filtrations $\F^n \subset \F$ in Theorem \ref{thlimsup1} and when we have the convergence of filtrations $\F^n \xrightarrow{w} \F$ in 
Theorem \ref{thlimsup2}. Finally, Theorem \ref{cvGamma} is proved.


\section{Applications}
\label{appl}
 
\subsection{Application to discretizations} 
Let us apply what we have proved in the case of discretizations.
\begin{prop}
\label{cvGammadiscr}
Let us consider a c\`adl\`ag process X such that \mbox{$\P[\Delta X_t \not= 0]=0$} for every t. 
Let \mbox{$(\pi^n=\{t_1^n, \ldots t_{k^n}^n\})_n$} be an increasing sequence of subdivisions of $[0,T]$ with mesh going to 0
$(|\pi^n| \xrightarrow[n \to +\infty]{} 0)$. 
We define the sequence of discretized processes $(X^n)_n$ by $\forall n$, $\forall t$,
\mbox{$X_t^n = \sum_{i=1}^{k^n-1} X_{t_i^n} 1_{t_i^n \leqslant t < t_{i+1}^n}$.}  \\
Then \mbox{$\Gamma_n(L) \xrightarrow[n \to +\infty]{} \Gamma(L)$.}
\end{prop}

\demo
Let us consider $\F^X$ the natural filtration for $X$, $\F$ the right continuous associated filtration and $(\F^n)_n$ the natural filtrations
for the $(X^n)_n$. \\
- $X^n \xrightarrow[n \to +\infty]{} X$ $a.s.$ then in probability. \\
- $\forall n, \F^n \subset \F^X \subset \F$ by definition of $X^n$. \\
- $\P[\Delta X_L \not= 0]=0$ by hypothesis. \\
- Using Proposition \ref{hypA}, Aldous' criterion is filled. \\
Then using Theorem \ref{cvGamma}, $\Gamma_n(L) \xrightarrow[n \to +\infty]{} \Gamma(L)$.
\findem

\subsection{Application to financial models}

We are going to apply the previous results to financial models. For a study about those models, see for example the book \citep{LL}. \\
 
We wish to find the price of an American call option at the best time of exercise for the buyer. We denote by $T$ the maturity date of this
call option. The market is composed of an asset with risk of price $S_t$ at time $t$ and an asset without risk of price $S^0_t$ at time
$t$. We assume that $S_t$ follows the stochastic differential equation $dS_t=S_t(\mu dt + \sigma dB_t)$ where $\mu$ and $\sigma$ are positive
reals and $(B_t)$ is a standart brownian motion. We also assume that $S_t^0$ is solution of the ordinary differential equation 
$dS_t^0=rS_t^0 dt$ where $r>0$. \\

We define the actualized price of the asset with risk by $\tilde{S_t}=e^{-rt}S_t$. Then, we have  
$d\tilde{S_t}=\tilde{S_t}(\lambda dt + \sigma dB_t)$ where $\lambda=\mu - r$. The solution of this equation is well known : 
$\tilde{S_t}=\tilde{S_0} exp(\lambda t - \sigma^2 t /2 + \sigma B_t).$ \\

The natural filtration for $\tilde{S}$ is the brownian filtration, denoted by $\F$. At the optimal exercice date, the price of the option is
given by the following value in optimal stopping of horizon $T$ for $\tilde{S}$ :
$$\Gamma^{\tilde{S}}(T)=\underset{\tau \in \mathcal{T}}{\sup} \E[\tilde{S_\tau}],$$
where $\mathcal{T}$ is the set of $\F$ stopping times bounded by $T$. \\

It is usual to approximate the model of Black and Scholes by a sequence of models of Cox-Ross-Rubinstein. \\

On an adapted space, we consider a sequence $(X_i)$ of independent Bernoulli variables such that $\forall i, \P[X_i=1]=\P[X_i=-1]=1/2$. 
For every $n \in \N^*$, we consider $B_{kT/n}^n=\sqrt{T/n} \sum_{i=1}^k X_i, k=0, \ldots, n$. We assume that
the actualized prices $\tilde{S}_{kT/n}^n$ of the asset with risk at time $kT/n$ are given by the linear equation
\mbox{$\Delta \tilde{S}_{(k+1)T/n}^n = \tilde{S}_{kT/n}^n (\lambda_n T/n + \sigma_n \Delta B_{(k+1)T/n}^n)$}
where $\Delta \tilde{S}_{(k+1)T/n}^n = \tilde{S}_{(k+1)T/n}^n - \tilde{S}_{kT/n}^n$ and 
$\Delta B_{(k+1)T/n}^n = B_{(k+1)T/n}^n - B_{kT/n}^n$.\\ 

We extend processes $B^n$ and $\tilde{S}^n$ to $[0,T]$ by the following : $B^n_t = B_{kT/n}^n$ if $kT/n \leqslant t < (k+1)T/n$ and 
$\tilde{S}^n_t = \tilde{S}_{kT/n}^n$ if $kT/n \leqslant t < (k+1)T/n$. \\

The natural filtration for $\tilde{S}^n$ is $\F^n$ such that \mbox{$\F^n_t = \sigma(B_{kT/n}^n, kT/n \leqslant t)$,} for every $t$. 
At the optimal exercise date, the value of the option is given by the following reduite of horizon $T$ associated to $\tilde{S}^n$ :
$$\Gamma^{\tilde{S}^n}(T)=\underset{\tau \in \mathcal{T}^n}{\sup} \E[\tilde{S}^n_\tau],$$
where $\mathcal{T}^n$ is the set of $\F^n$ stopping times bounded by $T$. \\

We assume that $\lambda_n \xrightarrow[n \to +\infty]{} \lambda$ and $\sigma_n \xrightarrow[n \to +\infty]{} \sigma$. \\

Using Donsker's Theorem, we have :
$$(B^n,\tilde{S}^n) \cvl (B,\tilde{S}).$$

According to the Skorokhod representation Theorem, we can find processes $(X,Y)$ and $((X^n,Y^n))_n$ such that $\forall n$, 
$(X^n,Y^n) \sim (B^n,\tilde{S}^n)$, $(X,Y) \sim (B, \tilde{S})$ and $(X^n,Y^n) \cvps (X, Y)$. \\

But, $\tilde{S}^n$ is a continuous function of $B^n$ and $(X^n,Y^n) \sim (B^n,\tilde{S}^n)$ so $Y^n$ is a continuous function of $X^n$. 
Hence, $Y^n$ and $X^n$ have the same natural filtration $\F^{X^n}=\F^{Y^n}$. Similarly, $X$ and $Y$ have the same natural filtration 
$\F^X=\F^Y$.\\

Moreover, $B$ is a process with independent increments, so also is $X$. Then, using \citep[Theorem 2]{cvfiltration}, as $X^n \cvp X$, we have 
the corresponding convergence of filtrations : $\F^{X^n} \xrightarrow{w} \F^X$. Hence, $\F^{Y^n} \xrightarrow{w} \F^Y$.\\

$Y$ and $\tilde{S}$ have the same law so $Y$ is quasi-left continuous. $Y^n \cvp Y$, $\F^{Y^n} \xrightarrow{w} \F^Y$ and $Y$ is 
quasi-left continuous, so using Theorem \ref{cvGamma}, we have
$$\Gamma^{Y^n}(T) \xrightarrow[n \to +\infty]{} \Gamma^{Y}(T)$$ where 
$\Gamma^{Y^n}(T)=\underset{\tau \in \mathcal{T}^{Y^n}}{\sup} \E[Y^n_\tau]$ with $\mathcal{T}^{Y^n}$ the set of $\F^{Y^n}$ stopping times
bounded by $T$ and $\Gamma^{Y}(T)=\underset{\tau \in \mathcal{T}^{Y}}{\sup} \E[Y_\tau]$ with $\mathcal{T}^{Y}$ the set of $\F^{Y}$ stopping
times bounded by $T$. \\

But according to Remark \ref{GammaX}, the value in optimal stopping only depends on the law of the process. Here, $Y$ and $\tilde{S}$ have 
the same law so $\Gamma^{Y}(T)=\Gamma^{\tilde{S}}(T)$ and $Y^n$ and $\tilde{S}^n$ have the same law so 
$\Gamma^{Y^n}(T)=\Gamma^{\tilde{S}^n}(T)$. Then, the sequence of values in optimal stopping associated to the models of Cox-Ross-Rubinstein 
converges to the value in optimal stopping of the model of Black and Scholes~:
$$\Gamma^{\tilde{S}^n}(T) \xrightarrow[n \to +\infty]{} \Gamma^{\tilde{S}}(T).$$



\begin{thebibliography}{10}

\bibitem[Aldous(1978)]{Aldous_Stop1}
Aldous, D., 1978.
\newblock Stopping times and tightness.
\newblock {Ann. Proba.}, 6(2):335--340.

\bibitem[Aldous(1989)]{Aldous_Stop2}
Aldous, D., 1989.
\newblock Stopping times and tightness. II.
\newblock {Ann. Proba.}, 17(2):586--595.

\bibitem[Aldous(1981)]{preprintAldous}
Aldous, D., 1981.
\newblock Weak convergence of stochastic processes for processes viewed in the
  Strasbourg manner.
\newblock Preprint, Statis. Laboratory Univ. Cambridge.

\bibitem[Baxter and Chacon(1977)]{BC}
Baxter, J.R., Chacon, R.V., 1977.
\newblock Compactness of stopping times.
\newblock { Z. Wahrscheinlichkeitstheorie verw. Gebiete}, 40 (3):169--181.

\bibitem[Billingsley(1999)]{Bill}
Billingsley, P., 1999.
\newblock { Convergence of Probability Measures, Second Edition}.
\newblock Wiley and Sons, New York.

\bibitem[Br\'emaud and Yor(1978)]{BY}
Br\'emaud, P., Yor, M., 1978.
\newblock Changes of filtrations and of probability measures.
\newblock { Z. Wahrscheinlichkeitstheorie verw. Gebiete}, 45:269--295.

\bibitem[Coquet, M\'emin and S\l omi\'nski(2001)]{cvfiltration}
Coquet, F., M\'emin, J., S\l omi\'nski, L., 2001.
\newblock On weak convergence of filtrations.
\newblock { S\'eminaire de probabilit\'es XXXV, Lectures Notes in
  Mathematics, Springer Verlag, Berlin Heidelberg New York}, 1755:306--328.

\bibitem[Hoover(1991)]{Hoover}
Hoover, D.N. , 1991.
\newblock Convergence in distribution and Skorokhod convergence for the general
  theory of processes.
\newblock { Probab. Theory Related Fields}, 89(3):239--259.

\bibitem[Jacod and M\'emin(1981)]{JacodMemin}
Jacod, J., M\'emin, J., 1981.
\newblock Sur un type de convergence interm\'ediaire entre la convergence en
  loi et la convergence en probabilit\'e.
\newblock { S\'eminaire de Probabilit\'es, XV, Lectures Notes in
  Mathematics, Springer, Berlin}, 850:529--546.

\bibitem[Jacod and Shiryaev(2002)]{JS}
Jacod, J., Shiryaev, A.N., 2002.
\newblock { Limit Theorems for Stochastic Processes, Second Edition}.
\newblock Springer Verlag, Berlin Heidelberg New York.

\bibitem[Lamberton and Lapeyre(1997)]{LL}
Lamberton, D., Lapeyre, B., 1997.
\newblock { Introduction au calcul stochastique appliqu\'e \`a la finance,
  Seconde Edition}.
\newblock Ellipses Edition Marketing, Paris.

\bibitem[Lamberton and Pag\`es(1990)]{cvreduites}
Lamberton, D., Pag\`es, G., 1990.
\newblock Sur l'approximation des r\'eduites.
\newblock { Ann. Inst. Henri Poincar\'e}, 26(2):331--355.

\bibitem[Meyer(1978)]{Meyer}
Meyer, P.A., 1978.
\newblock Convergence faible et compacit\'e des temps d'arr\^et d'apr\`es
  {B}axter et {C}hacon.
\newblock { S\'eminaire de Probabilit\'es, XII, Lectures Notes in
  Mathematics, Springer, Berlin}, 649:411--423.

\bibitem[Renyi(1963)]{Renyi}
Renyi, A., 1963.
\newblock On stable sequences of events.
\newblock {Sankya, Ser A}, 25:293--302.


\end{thebibliography}
\end{document}